\documentclass[11pt]{article}
\usepackage{amssymb,amsmath,latexsym}
\usepackage{pdfsync}

\begin{document}

\begin{doublespace}

\def\1{{\bf 1}}
\def\ind{{\bf 1}}
\def\nn{\nonumber}

\def\sA {{\cal A}} \def\sB {{\cal B}} \def\sC {{\cal C}}
\def\sD {{\cal D}} \def\sE {{\cal E}} \def\sF {{\cal F}}
\def\sG {{\cal G}} \def\sH {{\cal H}} \def\sI {{\cal I}}
\def\sJ {{\cal J}} \def\sK {{\cal K}} \def\sL {{\cal L}}
\def\sM {{\cal M}} \def\sN {{\cal N}} \def\sO {{\cal O}}
\def\sP {{\cal P}} \def\sQ {{\cal Q}} \def\sR {{\cal R}}
\def\sS {{\cal S}} \def\sT {{\cal T}} \def\sU {{\cal U}}
\def\sV {{\cal V}} \def\sW {{\cal W}} \def\sX {{\cal X}}
\def\sY {{\cal Y}} \def\sZ {{\cal Z}}

\def\bA {{\mathbb A}} \def\bB {{\mathbb B}} \def\bC {{\mathbb C}}
\def\bD {{\mathbb D}} \def\bE {{\mathbb E}} \def\bF {{\mathbb F}}
\def\bG {{\mathbb G}} \def\bH {{\mathbb H}} \def\bI {{\mathbb I}}
\def\bJ {{\mathbb J}} \def\bK {{\mathbb K}} \def\bL {{\mathbb L}}
\def\bM {{\mathbb M}} \def\bN {{\mathbb N}} \def\bO {{\mathbb O}}
\def\bP {{\mathbb P}} \def\bQ {{\mathbb Q}} \def\bR {{\mathbb R}}
\def\bS {{\mathbb S}} \def\bT {{\mathbb T}} \def\bU {{\mathbb U}}
\def\bV {{\mathbb V}} \def\bW {{\mathbb W}} \def\bX {{\mathbb X}}
\def\bY {{\mathbb Y}} \def\bZ {{\mathbb Z}}
\def\R {{\mathbb R}} \def\RR {{\mathbb R}}
\def\n{{\bf n}}

\newcommand{\expr}[1]{\left( #1 \right)}
\newcommand{\cl}[1]{\overline{#1}}
\newtheorem{thm}{Theorem}[section]
\newtheorem{lemma}[thm]{Lemma}
\newtheorem{defn}[thm]{Definition}
\newtheorem{prop}[thm]{Proposition}
\newtheorem{corollary}[thm]{Corollary}
\newtheorem{remark}[thm]{Remark}
\newtheorem{example}[thm]{Example}
\numberwithin{equation}{section}
\def\ee{\varepsilon}
\def\qed{{\hfill $\Box$ \bigskip}}
\def\NN{{\cal N}}
\def\AA{{\cal A}}
\def\MM{{\cal M}}
\def\BB{{\cal B}}
\def\CC{{\cal C}}
\def\LL{{\cal L}}
\def\DD{{\cal D}}
\def\FF{{\cal F}}
\def\EE{{\cal E}}
\def\QQ{{\cal Q}}
\def\RR{{\mathbb R}}
\def\R{{\mathbb R}}
\def\L{{\bf L}}
\def\K{{\bf K}}
\def\S{{\bf S}}
\def\A{{\bf A}}
\def\E{{\mathbb E}}
\def\F{{\bf F}}
\def\P{{\mathbb P}}
\def\N{{\mathbb N}}
\def\eps{\varepsilon}
\def\wh{\widehat}
\def\wt{\widetilde}
\def\pf{\noindent{\bf Proof.} }
\def\pff{\noindent{\bf Proof} }
\def\beq{\begin{equation}}
\def\eeq{\end{equation}}
\def\bee{\begin{equation}}
\def\eee{\end{equation}}
\def\osc{\mathrm{Osc}}

\title{\Large \bf
Boundary Harnack principle and Martin boundary at infinity for subordinate Brownian motions }
\author{{\bf Panki Kim}\thanks{This work was supported by Basic Science Research Program through the National Research Foundation of Korea(NRF)
grant funded by the Korea government(MEST)
(2012-0000940).} \quad {\bf Renming Song}\thanks{Research supported in part by a grant from the Simons
Foundation (208236).}
\quad and \quad {\bf Zoran Vondra\v{c}ek}\thanks{Supported in part by the MZOS
grant 037-0372790-2801.}  }

\date{ }
\maketitle

\begin{abstract}
In this paper we study the Martin boundary of unbounded open sets at infinity  for a large class of
subordinate Brownian motions.
We first prove that, for such subordinate Brownian motions,
the uniform boundary Harnack principle at infinity holds for arbitrary unbounded open sets.
Then we introduce the notion of $\kappa$-fatness  at infinity for open sets and show that
the Martin boundary at infinity of any such open set
consists of exactly one point and that point is a minimal Martin boundary point.
\end{abstract}



%%%%%%%%%%%%%%%%%%%%%%%%%%%%%%%%%%%%%%%%%%%%  Introduction and main results  %%%%%%%%%%%%%%%%%%%%%%%%%%%%%%%%%%%%%%%%%%%%%%%%%%%%%%

\section{Introduction and main results}

The study of the boundary Harnack principle for non-local operators started in the
late 1990's with \cite{B} which proved that the boundary Harnack principle
holds for the fractional Laplacian (or equivalently the rotationally invariant stable process) in bounded Lipschitz domains.
This boundary Harnack principle was extended to arbitrary open sets in \cite{SW}.
The final word in the case of the rotationally invariant $\alpha$-stable process was given in \cite{BKK}
where the so called uniform boundary Harnack principle was proved in arbitrary open sets with a
constant not depending on the set itself. Subsequently, the boundary Harnack principle was
extended to more general symmetric L\'evy processes, more precisely to subordinate Brownian
motions with ever more weaker assumptions on the Laplace exponents of the subordinators,
see \cite{KSV1}, \cite{KSV6}, \cite{KSV7} and \cite{KM2}.
Recently in \cite{BKuK}, a boundary Harnack principle was established in the setting of jump processes in metric measure spaces.

Let us be more specific and state the (slightly stronger) assumptions under which
the boundary Harnack principle was proved in \cite{KSV7}.
Let $S=(S_t)_{t\ge 0}$ be a subordinator  (a nonnegative L\'evy process with $S_0=0$) with Laplace exponent $\phi$ and
$W=(W_t,\P_x)_{t\ge 0, x\in \R^d}$ be a Brownian motion in $\R^d$, $d\ge 1$,  independent of $S$ with
$$
\E_x\left[e^{i
\xi \cdot (W_t-W_0)
}\right]=e^{-t{|\xi|^2}} \quad \xi\in \R^d, t>0.
$$
The process $X=(X_t,\P_x)_{t\ge 0, x\in \R^d}$ defined by $X_t:=W(S_t)$ is called a subordinate Brownian motion.
It is a rotationally invariant L\'evy process in $\R^d$ with characteristic exponent $\phi(|\xi|^2)$ and
infinitesimal generator $-\phi(-\Delta)$. Here $\Delta$ denotes the Laplacian and $\phi(-\Delta)$ is defined through functional calculus.

The function $\phi$ is a Bernstein function having the representation
$$
\phi(\lambda)=a+b\lambda +\int_{(0,\infty)}(1-e^{-\lambda t})\, \mu(dt)
$$
where $a,b\ge 0$ and $\mu$ is the measure satisfying $\int_{(0,\infty)}(1\wedge t)\, \mu(dt)<\infty$,
called the L\'evy measure of $\phi$ (or $S$).
Recall that $\phi$ is a complete Bernstein function if the measure $\mu$ has a completely
monotone density.
For basic facts about complete Bernstein functions, see \cite{SSV}.

Let us introduce the following upper and lower
scaling conditions on $\phi$ at infinity:

\medskip
\noindent
{\bf (H1):}
There exist constants $0<\delta_1\le \delta_2 <1$ and $a_1, a_2>0$  such that
\begin{equation}\label{e:H1}
a_1\lambda^{\delta_1} \phi(t) \le \phi(\lambda t) \le a_2 \lambda^{\delta_2} \phi(t), \quad \lambda \ge 1, t \ge 1\, .
\end{equation}
This is a condition on the asymptotic behavior of $\phi$ at infinity and it governs the behavior of the
subordinator $S$ for small time and small space (see \cite{KSV7, KSV8}).
Note that it follows from the second inequality above that $\phi$ has no drift, i.e., $b=0$.
Suppose that $\phi$ is a complete Bernstein function with the killing term
$a=0$ and that {\bf (H1)} holds.
The following boundary Harnack principle is proved in \cite[Theorem 1.1]{KSV7}:
There exists a constant $c=c(\phi,d)>0$ such that for every $
z \in \R^d$, every open set $D\subset \R^d$, every $r\in (0,1)$ and any nonnegative
functions $u, v$ on $\R^d$ which are regular harmonic in $D\cap B(
z, r)$ with respect to $X$ and vanish in $D^c \cap B(
z, r)$,
$$
\frac{u(x)}{v(x)}\,\le c\,\frac{u(y)}{v(y)}, \qquad x, y\in D\cap B(z, r/2).
$$
Here, and in the sequel,  $B(z,r)$ denotes
the open ball in $\R^d$
centered at $z$ with radius $r$.
This result was obtained as a simple consequence of the following approximate factorization
of nonnegative harmonic functions, see \cite[Lemma 5.5]{KSV7}: There exists a constant
$c=c(\phi, d)>1$ such that for every $z \in \R^d$, every open set $D \subset B(z,r)$ and any nonnegative function $u$ on $\RR^d$ which is regular harmonic in $D$ with respect to $X$ and
vanishes a.e.~in $D^c \cap B(z, r)$,
\begin{eqnarray}\label{e:approx-factor-finite}
c^{-1}\, \E_x[\tau_{D}] \int_{B(
z, r/2)^c}  j(|y-
z|)  u(y)\, dy \le u(x) \le
c \ \E_x[\tau_{D}] \int_{B(
z, r/2)^c} j(|y-
z|)  u(y)dy
\end{eqnarray}
for every  $ x \in D \cap B(
z, r/2)$. Here $j$ denotes the density of the L\'evy measure of $X$ and $\tau_D$ the first exit time of $X$ from $D$.
In the case of the rotationally invariant $\alpha$-stable process, \eqref{e:approx-factor-finite} is proved earlier in \cite{BKK}.

Note that the boundary Harnack principle is a result about the decay of non-negative harmonic
functions near the (finite) boundary points. It is an interesting problem to study the decay of
non-negative harmonic functions at infinity (which may be regarded as a ``boundary point at infinity''
of unbounded sets). This is the main topic of the current paper. In order to study the behavior of
harmonic functions at infinity, one needs large space and large time properties of the underlying
process $X$. This requires a different type of assumptions than {\bf (H1)} which gives only small
space and small time properties of $X$.
Therefore, in addition to {\bf (H1)}, we will assume the corresponding upper and lower scaling
conditions of $\phi$ near zero:

\noindent
{\bf (H2):}
There exist constants $0<\delta_3\le \delta_4 <1$ and $a_3, a_4>0$  such that
\begin{equation}\label{e:H2}
a_3\lambda^{\delta_4} \phi(t) \le \phi(\lambda t) \le a_4 \lambda^{\delta_3} \phi(t), \quad \lambda \le 1, t \le 1\, .
\end{equation}
This is a condition on the asymptotic behavior of $\phi$ at zero and it governs the
behavior of the subordinator $S$ for large time and large space
(see \cite{KSV8} for details and examples).

In the recent preprint \cite{KSV8} we studied the potential theory of subordinate
Brownian motions under the assumption that $\phi$ is a complete Bernstein function
satisfying both conditions {\bf (H1)} and {\bf (H2)}. We were able to extend many
potential-theoretic results that were proved under {\bf (H1)} (or similar assumptions
on the small time and small space behavior)  for radii $r\in (0,1)$ to the case of all
$r>0$ (with a uniform constant \emph{not} depending on $r>0$). In particular, we proved
a uniform boundary Harnack principle with explicit decay rate (in open sets satisfying
the interior and the exterior ball conditions) which is valid for all $r>0$.
The current paper is a continuation of \cite{KSV8} and is based on the results of \cite{KSV8}.

For any open set $D$, we use $X^D$ to denote the subprocess of $X$ killed upon exiting $D$.
In case $D$ is a Greenian open set in $\R^d$ we will use $G_D(x,y)$ to denote the Green function of $X^D$.
For a Greenian open set $D\subset \R^d$, let
$$
K_D(x,y):=  \int_{D} G_D(x,z) j(|z-y|)\, dz, \qquad (x,y) \in D \times \overline{D}^c
$$
be the Poisson kernel of $X$ in $D \times \overline{D}^c$.

The first goal of this paper is to prove the following  approximate factorization of regular harmonic functions vanishing at infinity.

\begin{thm}\label{bhp-inf} Suppose that $\phi$ is a complete Bernstein function satisfying
{\bf (H1)}--{\bf (H2)}, let $d >2(\delta_2 \vee \delta_4)$, and let $X$ be a rotationally
invariant L\'evy process in $\R^d$ with
characteristic exponent $\phi(|\xi|^2)$.  For every $a>1$, there exists $C_1=C_1(\phi, a)>1$
such that for any $r\ge 1$, any open set $U\subset \overline{B}(0,r)^c$ and any nonnegative
function $u$ on $\R^d$ which is regular harmonic with respect to $X$ in $U$ and vanishes a.e.
on $\overline{B}(0,r)^c\setminus U$, it holds that
\begin{equation}\label{e:bhp-inf}
C_1^{-1}K_U(x,0) \int_{B(0,ar)}u(z)\, dz \le u(x) \le C_1 K_U(x,0) \int_{B(0,ar)}u(z)\, dz \, ,
\end{equation}
for all $x\in U\cap \overline{B}(0,ar)^c$.
\end{thm}
Note that $K_U(x,0)=\int_U G_U(x,y) j(|y|)\, dy$.
A consequence of the assumption $d >2(\delta_2 \vee \delta_4)$ (always true for $d\ge 2$) in Theorem \ref{bhp-inf} is that the process $X$ is transient and points are polar.
Under this assumption, the Green function $G(x,y)$ of the process $X$ exists, and by \eqref{e:G} below we have  $G_U(x,y) \le G(x,y)
\asymp |x-y|^{-d}\phi(|x-y|^{-2})^{-1}$. This
will be used several times in this paper.

In the case of the rotationally invariant $\alpha$-stable process, Theorem \ref{bhp-inf}
(for $a=2$) was obtained in \cite[Corollary 3]{Kw} from \eqref{e:approx-factor-finite}
by using the inversion with respect to spheres and the Kelvin transform of harmonic functions
for the stable process. Since the Kelvin transform method works only for stable
processes we had to use a different approach to prove \eqref{e:bhp-inf}.
We followed the method used in \cite{KSV7} to prove \eqref{e:approx-factor-finite}, making
necessary changes at each step. The main technical difficulty of the proof is the
delicate upper estimate of the Poisson kernel $K_{\overline{B}(0,r)^c}(x,0)$ of
the complement of the ball given in Lemma \ref{l:poisson-kernel-estimate}, where
the full power of the results from \cite{KSV8} was used.

Theorem \ref{bhp-inf} gives  the following scale invariant boundary Harnack inequality at infinity.

\begin{corollary}[Boundary Harnack Principle at Infinity]\label{c:bhp-inf}
Suppose that $\phi$ is a complete Bernstein function satisfying
{\bf (H1)}--{\bf (H2)}, $d >2(\delta_2 \vee \delta_4)$, and that
$X$ is a rotationally invariant L\'evy process in $\R^d$ with
characteristic exponent $\phi(|\xi|^2)$.
For each $a>1$ there exists $C_2=C_2(\phi,a)>1$ such that for
any $r\ge 1$, any open set $U\subset \overline{B}(0,r)^c$ and any nonnegative functions
$u$ and $v$ on $\R^d$ that are regular harmonic in $U$ with respect to $X$ and
vanish a.e. on $\overline{B}(0,r)^c\setminus U$, it holds that
\begin{equation}\label{e:bhp-inf-cor}
    C_2^{-1} \frac{u(y)}{v(y)} \le  \frac{u(x)}{v(x)}\le C_2 \frac{u(y)}{v(y)}\, ,\qquad \textrm{for all }x,y\in U\cap \overline{B}(0,ar)^c\, .
\end{equation}
\end{corollary}

The boundary Harnack principle is the main tool in identifying the (minimal) Martin boundary
(with respect to the process $X$) of an open set. Recall that for $\kappa\in (0,1/2]$, an
open set $D$ is said to be $\kappa$-fat open at $Q\in \partial D$, if there exists $R>0$
such that for each $r\in (0,R)$ there exists a point $A_r(Q)$ satisfying
$B(A_r(Q), \kappa r)\subset D\cap B(Q,r)$. If $D$ is $\kappa$-fat at each
boundary point $Q\in \partial D$ with the same $R>0$, $D$ is called
$\kappa$-fat with characteristics $(R,\kappa)$. In case $X$ is a subordinate
Brownian motion via a subordinator with a complete Bernstein Laplace exponent
regularly varying at infinity with index in $(0,1)$,
it is shown in \cite{KSV1} that the minimal Martin boundary of a bounded $\kappa$-fat
open set can be identified with the Euclidean boundary.

Corollary \ref{c:bhp-inf} enables us to identify the Martin boundary and the
minimal Martin boundary at infinity of a large class of  open sets  with respect to $X$.
To be more precise, let us first define $\kappa$-fatness at infinity.

\begin{defn}\label{d:kappa-fat-infty}
Let $\kappa \in (0,1/2]$. We say that an open set $D$ in $\R^d$ is $\kappa$-fat at infinity
if there exists $R>0$ such that for every $r\in [R,\infty)$ there exists $A_r \in \R^d$ such that $B(A_r, \kappa r)\subset D\cap \overline{B}(0,r)^c$ and $|A_r|< \kappa^{-1} r$. The pair $(R,\kappa)$ will be called the characteristics of the $\kappa$-fat open set $D$ at infinity.
\end{defn}

Note that all half-space-like open sets, all exterior open sets and all infinite cones are $\kappa$-fat at infinity.
Examples of  disconnected open sets which are $\kappa$-fat at infinity are
\begin{description}
\item{(i)}    $\{x=(x_1,\dots, x_{d-1},x_d)\in \R^d: x_d<0  \text{ or } x_d>1\}$;
\item{(ii)}  $
\bigcup^\infty_{n=1}B(x^{(n)}, 2^{n-2})$ with $|x^{(n)}|=2^n$.
\end{description}

Let $D\subset \R^d$ be an open set which is $\kappa$-fat at infinity with characteristics $(R,\kappa)$. Fix $x_0\in D$ and define
$$
M_D(x, y):=\frac{G_D(x, y)}{G_D(x_0, y)}, \qquad x, y\in D,~y\neq x_0.
$$
As the process $X^D$ satisfies Hypothesis (B) in \cite{KW}, $D$  has a Martin boundary $\partial_M D$
with respect to $X$ and $M_D(x ,\, \cdot\,)$ is  continuously extended  to $\partial_M D$.
A point $w \in \partial_M D$ is called an infinite  Martin boundary point if every sequence $(y_n)_{n\ge 1}$, $y_n\in D$, converging to $w$ in the Martin topology is unbounded (in the Euclidean metric). The set of all infinite Martin boundary points will be denoted by $\partial_M^{\infty} D$ and we call this set the Martin boundary at infinity.

By using the boundary Harnack principle at infinity we first show that if $D$ is $\kappa$-fat at infinity, then there exists the limit
\begin{equation}\label{e:martin-kernel-limit}
M_D(x,\infty)=\lim_{y\in D, |y|\to \infty} M_D(x,y)\, .
\end{equation}
The existence of this limit shows that $\partial_M^{\infty} D$ consists of
a single point which we denote by $\partial_{\infty}$. Finally, we prove that
$\partial_{\infty}$ is a minimal Martin boundary point. These findings are summarized
in the second main result of the paper.

\begin{thm}\label{mM-inf}
Suppose that $\phi$ is a complete Bernstein function satisfying {\bf (H1)}--{\bf (H2)}, $d >2(\delta_2 \vee \delta_4)$, and $X$ is a rotationally invariant L\'evy process in $\R^d$ with
characteristic exponent $\phi(|\xi|^2)$.
Then the Martin boundary at infinity with respect to $X$ of any open set $D$
which is $\kappa$-fat at infinity
consists of exactly one point $\partial_{\infty}$. This point is a minimal Martin boundary point.
\end{thm}

We emphasize that this result is proved without any assumption on the finite boundary points. In particular, we do not assume that $D$ is $\kappa$-fat.

To the best of our knowledge, the only case where the Martin boundary at infinity has been
identified is the case of the rotationally invariant $\alpha$-stable process, see \cite{BKK}.
Again, the Kelvin transform method was used to transfer results for finite
boundary points to the infinite boundary point.
As we have already pointed out,
the Kelvin transform is not available for more general processes.

We remark here that for one-dimensional L\'evy processes (satisfying much weaker assumptions than ours) it is  proved in
\cite[Theorem 4]{Sil}  that the minimal Martin boundary at infinity for the half-line $D=(0,\infty)$ is one point. The question of the Martin boundary at infinity is not addressed in \cite{Sil}.

The paper is organized as follows. In the next section we introduce necessary notation and definitions,
and recall some results that follow from {\bf (H1)} and
{\bf (H2)} obtained in \cite{KSV8}.
Section 3 is devoted to the proofs of Theorem \ref{bhp-inf} and Corollary \ref{c:bhp-inf}. At the end of the section we collect some consequences of these two results. In the first part of Section 4 we study non-negative harmonic functions in unbounded sets that are $\kappa$-fat at infinity. The main technical result is
the oscillation reduction in
Lemma \ref{l:oscillation-reduction} immediately leading to \eqref{e:martin-kernel-limit}. Next we look at
the Martin and the minimal Martin boundary at infinity and give a proof of Theorem \ref{mM-inf}. We finish the paper by discussing the Martin boundary of the half-space.

Throughout this paper, the constants $C_1$, $C_2$,  $C_3, \dots$ will be fixed. The lowercase constants $
c_1, c_2, \dots$ will denote generic constants whose exact values are not important and can change from one appearance to another.
The dependence of the lower case constants on the dimension $d$ and the function $\phi$ may not be mentioned explicitly.
The constant $c$ that depends on the parameters $\delta_i$ and $a_i$, $i=1,2,3,4$, appearing in {\bf (H1)} and {\bf (H2)} will be simply denoted as $c=c(\phi)$.
We will use ``$:=$" to denote a definition, which is read as ``is defined to
be". For $a, b\in \bR$, $a\wedge b:=\min \{a, b\}$ and $a\vee
b:=\max\{a, b\}$.

For any open set $U$, we denote by $\delta_U (x)$ the
distance between $x$ and the complement of $U$, i.e.,
$\delta_U(x)=\text{dist} (x, U^c)$.
For functions $f$ and $g$, the notation ``$f\asymp g$"
means that there exist constants $c_2
 \ge c_1>0$ such that
$c_1 \, g \leq f \leq c_2 \, g$.
For every
function $f$, we extend its definition to the cemetery point $\partial$ by setting
$f(\partial )=0$.
For every function $f$, let $f^+:=f \vee 0$.
We will use $dx$ to denote the
Lebesgue measure in $\bR^d$ and, for a Borel set $A\subset \bR^d$, we
also use $|A|$ to denote its Lebesgue measure.
We denote $\overline{B}(0,r)^c:=\{y \in \R^d: |x-y| >r\}$.
Finally, for a point $x=(x_1,\dots, x_{d-1}, x_d)\in \R^d$ we sometimes write $x=(\wt{x},x_d)$ with $\wt{x}=(x_1,\dots, x_{d-1})\in \R^{d-1}$.

%%%%%%%%%%%%%%%%%%%%%%%%%%%%%%%%%%%%%%%%%%%%%%%%%%%%%%%  Preliminaries  %%%%%%%%%%%%%%%%%%%%%%%%%%%%%%%%%%%%%%%%%%%%%%%%%%%%%%%%%%

\section{Preliminaries}

In this section we recall  some results from \cite{KSV8}.
Recall that a function $\phi:(0,\infty)\to (0,\infty)$ is a Bernstein function if
it is $C^{\infty}$ function on $(0,\infty)$ and $(-1)^{n-1} \phi^{(n)}\ge 0$ for all $n\ge 1$.
It is well known that, if $\phi$ is a Bernstein function, then
\begin{equation}\label{e:Berall}
\phi(\lambda t)\le \lambda\phi(t) \qquad \text{ for all } \lambda \ge 1, t >0\, .
\end{equation}

Clearly \eqref{e:Berall} implies the following observation.
\begin{lemma}\label{l:phi-property}
If $\phi$ is a Bernstein function, then every $\lambda >0$,
$$
1 \wedge  \lambda\le \frac{\phi(\lambda t)}{\phi(t)} \le 1 \vee \lambda\, ,\quad \textrm{for all }t>0\, .
$$
\end{lemma}

 Note that with this lemma, we can replace expressions of the type $\phi(\lambda t)$, when $\phi$ is a Bernstein function,
with $\lambda >0$ fixed and $t>0$ arbitrary, with $\phi(t)$ up to a multiplicative
constant depending on $\lambda$. We will often do
this without explicitly mentioning it.

Recall that a subordinator $S=(S_t)_{t\ge 0}$ is simply a nonnegative L\'evy process
with $S_0=0$.
Let $S=(S_t)_{t\ge 0}$ be a subordinator with Laplace exponent $\phi$.
The function $\phi$ is a
Bernstein function with $\phi(0)=0$ so it has
the representation
$$
\phi(\lambda)= b\lambda +\int_{(0,\infty)} (1-e^{-\lambda t})\, \mu(dt)\, ,
$$
where $b\ge 0$ is the drift and $\mu$ the L\'evy measure of $S$.

A Bernstein function $\phi$ is a complete Bernstein function if its
L\'evy measure $\mu$ has a completely monotone density, which will
be denoted by $\mu(t)$. Throughout this paper we assume that $\phi$ is a complete Bernstein function.
In this case, the potential measure $U$ of $S$ admits a completely monotone density $u(t)$ (cf.~\cite{SSV}).

Conditions {\bf (H1)}--{\bf (H2)} imply that
\begin{equation}\label{e:sc1}
c^{-1} \left(\frac{R}{r}\right)^{\delta_1 \wedge \delta_3} \le \frac{\phi(R)}{\phi(r)} \le c\left(\frac{R}{r}\right)^{\delta_2 \vee \delta_4}, \quad 0<r<R<\infty\, .
\end{equation}
(See \cite{KSV8} for details.)
Using \eqref{e:sc1}, we have the following  result which will be used many times later in the paper.
(See the proof of \cite[Lemma 4.1]{KSV7} for similar computations.)
\begin{lemma}[\cite{KSV8}]\label{l:integral-estimates-phi}
    Assume {\bf (H1)} and {\bf (H2)}.
    There exists a constant
        $c=c(\phi) \ge 1$
     such that
    \begin{eqnarray}
    \int_0^{\lambda^{-1}} \phi (r^{-2})^{1/2}\, dr    \le  c  \lambda^{-1}\phi(\lambda^{2})^{1/2},  &\quad &\textrm{for all }\lambda > 0\, ,\label{e:ie-1} \\
    \lambda^2 \int_0^{\lambda^{-1}} r \phi (r^{-2})\, dr + \int_{\lambda^{-1}}^\infty r^{-1} \phi (r^{-2})\, dr \le  c \phi(\lambda^2)\, ,&\quad &\textrm{for all } \lambda >0\, ,\label{e:ie-2}
    \\
  c^{-1} \phi(\lambda^2)^{-1} \le  \int_0^{\lambda^{-1}} r^{-1}\phi(r^{-2})^{-1}\, dr \le c \phi(\lambda^2)^{-1}\, , &\quad &\textrm{for all } \lambda >0\,.\label{e:ie-3}
    \end{eqnarray}
\end{lemma}

Recall that $S=(S_t)_{t\ge 0}$ is a subordinator with Laplace exponent $\phi$.
Let $W=(W_t)_{t\ge 0}$ be a $d$-dimensional Brownian motion,
$d\ge 1$, independent of $S$ and with transition density
$$
q(t,x,y)=(4\pi t)^{-d/2} e^{-\frac{|x-y|^2}{4t}}\, ,\quad x,y\in \R^d, \ t>0\, .
$$
The process $X=(X_t)_{t\ge 0}$ defined by $X_t:=W(S_t)$ is called a subordinate Brownian motion.
$X$ is a rotationally invariant L\'evy process with characteristic exponent $\phi(|\xi|^2)$, $\xi\in \R^d$.
Throughout this paper $X$ is always such a subordinate Brownian motion.
The L\'evy measure of $X$ has a density $J(x)=j(|x|)$ where $j:(0,\infty)\to (0,\infty)$ is given by
$$
j(r):=\int_0^{\infty}(4\pi t)^{-d/2} e^{-r^2/(4t)}\mu(t)\, dt\, .
$$
Note that $j$ is continuous and decreasing.
Recall that the infinitesimal generator $\sL$ of the process $X$ (e.g. \cite[Theorem 31.5]{Sat}) is given by
\begin{equation}\label{e:infinitesimal-generator}
    \sL f(x)=\int_{\R^d}\left( f(x+y)-f(x)-y\cdot \nabla f(x) \1_{\{|y|\le \eps \}} \right)\, J(y)dy
\end{equation}
for every $\eps>0$ and $f\in C_b^c(\R^d)$,
where $C_b^c(\R^d)$ is the collection of $C^2$ functions which,
along with its partial derivatives of up to order 2, are bounded.

By the Chung-Fuchs criterion the process $X$ is transient if and only if
$$
\int_0^1 \frac{\lambda^{d/2-1}}{\phi(\lambda)}\, d\lambda <\infty\, .
$$
It follows that $X$ is always transient when $d\ge 3$.
In case {\bf (H2)} holds, $X$ is transient provided $\delta_4 < d/2$ (which is true if $d\ge 2$).
When $X$ is transient the occupation measure of $X$ admits a density $G(x,y)$
which is called the Green function of $X$ and is given by the formula
$G(x,y)=g(|x-y|)$ where
\begin{equation}\label{e:green-function}
g(r):=\int_0^{\infty}(4\pi t)^{-d/2} e^{-r^2/(4t)}u(t)\, dt\, .
\end{equation}
Here $u$ is the potential density of the subordinator $S$.
Note that by the transience assumption, the integral converges.
Moreover, $g$ is continuous and decreasing.
Furthermore, {\bf (H1)}--{\bf (H2)} imply the following estimates.
\begin{thm}[\cite{KSV8}]\label{t:J-G}
Assume  both {\bf (H1)} and {\bf (H2)}.
\begin{itemize}
    \item[(a)] It holds that
    \begin{equation}\label{e:J}
    j(r)\asymp r^{-d}\phi(r^{-2})\, , \quad \textrm{for all }r> 0\, .
    \end{equation}
    \item[(b)] If $d > 2(\delta_2\vee \delta_4)$ then the process $X$ is transient  and it holds
    \begin{equation}\label{e:G}
    g(r)\asymp r^{-d}\phi(r^{-2})^{-1}\, , \quad \textrm{for all }r> 0\, .
    \end{equation}
\end{itemize}
\end{thm}

As a consequence of \eqref{e:J}, we have
\begin{corollary}\label{c:doubling-condition}
Assume {\bf (H1)} and {\bf (H2)}. For every $L>1$, there exists a constant $c=c(L)>0$ such that
\begin{equation}\label{e:doubling-condition}
j(r)\le c j(L r)\, ,\quad r>0\, .
\end{equation}
\end{corollary}

For any open set $D$, we use $\tau_D$ to denote the first exit
time of $D$, i.e., $\tau_D=\inf\{t>0: \, X_t\notin D\}$.
Given  an open set $D\subset \R^d$, we define
$X^D_t(\omega)=X_t(\omega)$ if $t< \tau_D(\omega)$ and
$X^D_t(\omega)=\partial$ if $t\geq  \tau_D(\omega)$, where
$\partial$ is a cemetery state.

Let $p(t, x, y)$ be the transition density of $X$.
By the strong Markov property,
$$
p_D(t, x, y)\,:=\,p(t, x, y)\,-\,\E_x[ \,p(t-\tau_D, X_{\tau_D}, y)\,;\, t>\tau_D]\, , \qquad x, y \in D\, ,
$$
is the transition density of $X^D$.
A subset $D$ of $\R^d$ is said to be Greenian (for $X$) if $X^{D}$ is transient.
For a Greenian set $D\subset \R^d$, let $G_D(x,y)$ denote the Green function of $X^D$,
i.e.,
$$
G_D(x, y):=\int^\infty_0p_D(t, x, y)dt \qquad \text{for }\,\,x, y\in D.
$$

We define the Poisson kernel $K_D(x,z)$ of $D$ with respect to $X$ by
\begin{equation}\label{e:poisson-kernel-def}
K_D(x,z)=\int_{\overline{D}^c} G_D(x,y)J(y,z)\, dy,  \quad (x,z)\in D \times
\overline{D}^c \, .
\end{equation}
Then by \cite[Theorem 1]{IW} we get that for every
Greenian open subset $D$, every nonnegative Borel measurable function $f \ge 0$ and $x \in D$,
\begin{equation}\label{newls}
\E_x\left[f(X_{\tau_D});\,X_{\tau_D-} \not= X_{\tau_D}  \right]
=\int_{\overline{D}^c} K_D(x,y)f(y)dy.
\end{equation}
Using the continuities of $G_D$ and $J$, one can easily check
that $K_D$ is continuous on $D \times
\overline{D}^c$.

Equations \eqref{e:J} and \eqref{e:poisson-kernel-def}
 give the following estimates on
the Poisson kernel of $B(x_0,r)$ for all $r>0$.
\begin{prop}
[\cite{KSV8}]
\label{p:poisson-kernel-estimate}
Assume {\bf (H1)} and {\bf (H2)}.
There exist $c_1=c_1 (\phi)>0$ and $c_2=c_2 (\phi)>0$ such that for every $r >0$ and $x_0 \in \R^d$,
\begin{eqnarray}
K_{B(x_0,r)}(x,y) \,&\le &\, c_1 \, j(|y-x_0|-r) \left(\phi(r^{-2})\phi((r-|x-x_0|)^{-2})\right)^{-1/2} \label{e:pke-upper}\\
 &\le &\, c_1 \, j(|y-x_0|-r) \phi(r^{-2})^{-1}\ \label{e:pke-upper-worse}
\end{eqnarray}
for all $(x,y) \in B(x_0,r)\times \overline{B(x_0,r)}^c$ and
\begin{equation}\label{e:pke-lower}
K_{B(x_0, r)}(x_0, y) \,\ge\, c_2\, j(|y-x_0|) \phi(r^{-2})^{-1}, \qquad \textrm{ for all } y \in \overline{B(x_0, r)}^c.
\end{equation}
\end{prop}

To discuss the Harnack inequality and the boundary Harnack principle, we first recall the definition of harmonic functions.
\begin{defn}\label{D:1.1}
A function $f: \R^d\to [0, \infty)$ is said to be
\begin{description}
\item{(1)} harmonic in an open set $D\subset \R^d$ with respect to
$X$ if for every open set $B$ whose closure is a compact subset of $D$,
\begin{equation}\label{e:har}
f(x)= \E_x \left[ f(X(\tau_{B}))\right] \qquad
\mbox{for every } x\in B;
\end{equation}
\item{(2)} regular harmonic in $D$ for $X$ if
for each $x \in D$,
$f(x)= \E_x\left[f(X(\tau_{D}));\, \tau_D<\infty\right]$;
\item{(2)}
 harmonic for $X^D$ if it is harmonic for $X$ in $D$ and
vanishes outside $D$.
\end{description}
\end{defn}
We note that, by the strong Markov property of $X$,
every regular harmonic function is automatically harmonic.

Under the assumptions {\bf (H1)} and {\bf (H2)},
the following uniform Harnack inequality and uniform boundary Harnack principle
for all $r>0$ are true.

\begin{thm}[\cite{KSV8}]\label{t:uhi}
Assume {\bf (H1)} and {\bf (H2)}.
There exists $c=c(\phi)>0$ such that, for any $r>0$, $x_0\in \R^d$, and any function $u$ which is nonnegative on $\R^d$ and harmonic with respect to $X$ in $B(x_0, r)$, we have
$$
u(x)\,\le\, c\, u(y), \quad \textrm{for all }x, y\in B(x_0, r/2).
$$
\end{thm}
\begin{thm}[\cite{KSV8}]\label{t:ubhp}
Assume {\bf (H1)} and {\bf (H2)}.
There exists a constant $c= c(\phi
)>0$ such that for every $z_0 \in \R^d$, every open set $D\subset \R^d$, every $r >0$
and any nonnegative functions $u, v$ in $\R^d$ which are regular harmonic
in $D\cap B(z_0, r)$ with respect to $X$ and vanish in $D^c \cap B(z_0, r)$, we have
        $$
        \frac{u(x)}{v(x)}\,\le c\,\frac{u(y)}{v(y)}, \quad \mbox{ for all } x, y\in D\cap B(z_0, r/2).
        $$
\end{thm}

For $x\in \bR^d$, let $\delta_{\partial D}(x)$ denote the Euclidean
distance between $x$ and $\partial D$. Recall that  $\delta_{ D}(x)$ is the Euclidean
distance between $x$ and $D^c$.

We say that an open set $D$ satisfies both the {\it uniform interior ball condition} and the {\it uniform exterior ball condition} with
 radius $R$ for every $x\in D$
radius $R$ if for every $x\in D$
with $\delta_{\partial D}(x)\le R
$ and $y\in \bR^d \setminus
\overline D$ with $\delta_{\partial D}(y) \le R
$, there are $z_x,
z_y\in \partial D$ so that $|x-z_x|=\delta_{\partial D}(x)$,
$|y-z_y|=\delta_{\partial D}(y)$ and that $B(x_0, R
)\subset D$ and
$B(y_0, R
)\subset \bR^d \setminus \overline D$, where
$x_0=z_x+R
(x-z_x)/|x-z_x|$ and $y_0=z_y+R
(y-z_y)/|y-z_y|$.

The following is the one of main results in \cite{KSV8} --
the global  uniform
boundary Harnack principle with explicit decay rate  on open sets in $\bR^d$ with
the interior and exterior ball conditions with
radius $R$ for all $R >0$.

\begin{thm}[\cite{KSV8}]\label{L:2}
Assume {\bf (H1)} and {\bf (H2)}. There exists
$c=c(\phi)>0$  such that for every open set $D$ satisfying the interior and exterior ball conditions with
radius $R >0$,  $r \in (0, R]$, every $Q\in \partial D$ and every nonnegative function $u$ in $\R^d$
which is harmonic in $D \cap B(Q, r)$ with respect to $X$ and
vanishes continuously on $ D^c \cap B(Q, r)$, we have
\begin{equation}\label{e:bhp_m}
\frac{u(x)}{u(y)}\,\le c\,
\sqrt{\frac{\phi(\delta_{D}(y)
^{-2})}{\phi(\delta_{D}(x)
^{-2})}}
\qquad \hbox{for every } x, y\in  D \cap B(Q, \frac{r}{2}).
\end{equation}
\end{thm}

%%%%%%%%%%%%%%%%%%%%%%%%%%%%%%%%%%%%%%%%%%%%%%%  Boundary Harnack principle at infinity  %%%%%%%%%%%%%%%%%%%%%%%%%%%%%%%%%%%%%%%%%%

\section{Boundary Harnack principle at infinity}

The goal of this section is to prove the
scale invariant boundary Harnack principle at infinity (Theorem \ref{bhp-inf} and Corollary \ref{c:bhp-inf}).
In the remainder of this paper we assume that $\phi$
is a complete Bernstein function
satisfying  {\bf (H1)}--{\bf (H2)} and $d >2(\delta_2 \vee \delta_4)$, and
that $X$ is a rotationally invariant L\'evy process in $\R^d$ with
characteristic exponent $\phi(|\xi|^2)$.
Under these assumptions,
by \eqref{e:G},
$g$ satisfies the following property which we will use frequently: For every $L>1$, there exists
$c=c(L, \phi)>0$ such that
\begin{equation}\label{e:doubling-conditionG}
g(r)\,\le\, c\, g(Lr)\, ,\quad r>0\, .
\end{equation}

To prove Theorem \ref{bhp-inf} we need several lemmas.
For $x\in \R^d$ and $0<r_1<r_2$, we use $A(x,r_1,r_2)$ to denote the annulus
$\{y\in \R^d:\, r_1 <|y-x|  \le r_2\}$.

\begin{lemma}\label{l-exit}
For every $a \in (1,\infty)$, there exists
$c=c(\phi, a)>0$
such that for any $r>0$ and any open set $D\subset \overline{B}(0, r)^c$ we have
$$
{\P}_x\left(X_{\tau_D} \in \overline{B}(0, r)\right) \,\le\, c\, r^d K_D(x,0)\, , \qquad x \in D\cap {B}(0, ar)^c\, .
$$
\end{lemma}
\pf Let $\psi\in C^{\infty}_c(\R^d)$ be a function such that $0\le \psi\le 1$,
$$
\psi(y)=\left\{\begin{array}{ll}
    0\, ,& |y|>\frac{a+1}{2}\, ,\\
    1\, ,& |y|\le 1\, ,
    \end{array}\right.
$$
and $\sup_{y\in \R^d}
\sum_{i,j=1}^d
\left|\frac{\partial^2}{\partial y_i \partial y_j}\psi(y)\right|\le c_1=c_1(a)$. For $r>0$ define $\psi_r(y):=\psi(y/r)$. Then $\psi\in C^{\infty}_c(\R^d)$, $0\le \psi_r \le 1$,
$$
\psi_r(y)=\left\{\begin{array}{ll}
    0\, ,& |y|>\frac{a+1}{2}r\, ,\\
    1\, ,& |y| \le r\, ,
    \end{array}\right.
$$
and $\sup_{y\in \R^d}
\sum_{i,j=1}^d
\left|\frac{\partial^2}{\partial y_i \partial y_j}\psi_r(y)\right|\le c_1 r^{-2}$.

Let $x\in  D\cap {B}(0, ar)^c$. Recall that $\sL$ denotes the infinitesimal generator of $X$ and is given by \eqref{e:infinitesimal-generator}. Since $\psi_r(x)=0$ and $D\subset \overline{B}(0,r)^c$, by Dynkin's formula
(see, for instance, \cite[(5.8)]{Dyn})
we have
\begin{align}
   & \E_x\left[ \psi_r(X_{\tau_D})\right] = \int_D G_D(x,z)  {\sL}  \psi_r(z)\, dz\nonumber \\
    &=\int_{D\cap A(0,r, (a+2)r)} G_D(x,z)  {\sL}  \psi_r(z)\, dz + \int_{D\cap \overline{B}(0,(a+2)r)^c} G_D(x,z)  {\sL}  \psi_r(z)\, dz.
 \label{e:dynkin}
\end{align}
For $z\in D\cap A(0,r, (a+2)r)$ we
have
\begin{eqnarray*}
    |{\sL}\psi_r(z)|&=&\left|\int_{\R^d}\left(\psi_r(z+y)-\psi_r(z)-\nabla\psi_r(z)\cdot y \ind_{\{|y|\le r\}}\right)j(|y|)\, dy\right|\\
    &\le &\int_{\{|y|\le r\}} \left|\psi_r(z+y)-\psi_r(z)-\nabla\psi_r(z)\cdot y  \right| j(|y|)\, dy + 2\int_{\{r<|y|\}}j(|y|)\, dy\\
    &\le &\frac{c_2}{r^2}\int_{\{|y|\le r\}} |y|^2 j(|y|)\, dy + 2\int_{\{r<|y|\}}j(|y|)\, dy\\
    &\le & c_3\left( r^{-2} \int_0^r t\phi(t^{-2})\, dt + \int_r^{\infty} t^{-1} \phi(t^{-2})\, dt\right)\, ,
\end{eqnarray*}
where in the last line we have used \eqref{e:J}.
Thus by using \eqref{e:ie-2}, we get
$|{\sL}\psi_r(z)|\le c_4 \phi(r^{-2})$.
By Lemma \ref{l:phi-property} we see that $\phi(r^{-2})$ and $\phi((a+2)^{-2}r^{-2})$ are comparable (with a constant depending on $\phi$ and $a$). Therefore
\begin{equation}\label{e:small-z}
    |{\sL}\psi_r(z)|\le c_4\phi(|z|^{-2})\le c_5 |z|^d j(|z|)\le c_6 r^d j(|z|)\, , \quad r<|z|<(a+2)r.
\end{equation}

Now assume that $z\in D\cap \overline{B}(0,(a+2)r)^c$. Then $\psi_r(z)=0$ and $\nabla \psi_r(z)=0$ (note that $\psi_r$ is zero in a neighborhood of $z$). Therefore
\begin{eqnarray*}
    {\sL}\psi_r(z)&=& \int_{\R^d}\left(\psi_r(z+y)-\psi_r(z)-\nabla\psi_r(z)\cdot y \ind_{\{|y|\le r\}}\right)j(|y|)\, dy\\
    &=&\int_{\R^d}\psi_r(z+y)j(|y|)\, dy =\int_{|z+y|\le \frac{a+1}{2}r} \psi_r(z+y)j(|y|)\, dy\, ,
\end{eqnarray*}
where the last equality follows from the fact that $\psi_r(z+y)\neq 0$ only if $|z+y|\le \frac{a+1}{2}r$. Suppose that $|z+y|\le \frac{a+1}{2}r$. By the triangle inequality,
$$
    |y|\ge |z| -\frac{a+1}{2}r >|z|-\frac{a+1}{2}\, \frac{1}{a+2}|z|=\frac{a+3}{2(a+2)}|z|\ge \frac{|z|}{2}\, .
$$
It follows from  \eqref{e:doubling-condition} that $j(|y|)\le j(|z|/2)\le c_7 j(|z|)$. This implies that
\begin{equation}\label{e:large-z}
    {\sL}\psi_r(z) \le c_7 j(|z|) \int_{ |z+y|\le \frac{a+1}{2}r} dy \le c_8 r^d j(|z|)\, .
\end{equation}
Combining \eqref{e:dynkin}--\eqref{e:large-z}
we obtain
\begin{eqnarray*}
    \E_x\left[ \psi_r(X_{\tau_D})\right]
    &=&\int_{D\cap A(0,r, (a+2)r)} G_D(x,z)  {\sL}  \psi_r(z)\, dz + \int_{D\cap \overline{B}(0,(a+2)r)^c} G_D(x,z)  {\sL}  \psi_r(z)\, dz\\
    &\le &c_9 r^d \int_D G_D(x,z)j(|z|)\, dz
    = c_9 r^d K_D(x,0)\, .
\end{eqnarray*}
Finally, since $
{\bf 1}_{\overline{B(0,r)}}\le \psi_r$,
$
    \P_x(X_{\tau_D}\in \overline{B(0,r)})\le \E_x\left[ \psi_r(X_{\tau_D})\right] \le c_9 r^d K_D(x,0)\, .
$
\qed

\begin{lemma}\label{l:poisson-kernel-estimate}
Let $1<p<q<\infty $. There exists $c=c(\phi, p,q)>1$ such that for all
$r\ge 1/4$ it holds that
\begin{equation}\label{e:poisson-kernel-estimate}
     K_{\overline{B}(0,r)^c} (x,z)\le c
    r^{-d}\left(\phi(r^{-2})^{-1/2} \phi((r-|z|)^{-2})^{1/2}+1\right)
\end{equation}
for all $x\in A(0,pr, qr)$ and $z\in B(0,r)$.
\end{lemma}
\pf
We rewrite the Poisson kernel $K_{\overline{B}(0,r)^c}$ as follows:
\begin{align*}
   & K_{\overline{B}(0,r)^c} (x,z)\,=\,\int_{\overline{B}(0,r)^c} G_{\overline{B}(0,r)^c}(x,y)j(|y-z|)\, dy \\
    &=\left(\int_{A(0,r,2qr)\cap\{|x-z|\le 2|x-y|\}} + \int_{A(0,r,2qr)\cap\{|x-z|> 2|x-y|\}}+ \int_{\overline{B}(0,2qr)^c} \right)
    G_{\overline{B}(0,r)^c}(x,y)j(|y-z|)\, dy\\
    &=:\, I_1+I_2+I_3\, .
\end{align*}
Note that, since $x\in A(0,pr, qr)$ and $z\in B(0,r)$,
for  $y \in A(0,r,2qr)$ with $|x-z|\le 2|x-y|$, we have
\begin{equation} \label{e:newge}
|x-y| \ge \frac12(|x|-|z|) \ge \frac12(p-1)r >(4q)^{-1}(p-1)\delta_{\overline{B}(0,r)^c}(y).
\end{equation}
We claim that when $y \in A(0,r,2qr)$ satisfies $|x-z|\le 2|x-y|$,
\begin{equation} \label{e:newge2}
G_{\overline{B}(0,r)^c}(x,y) \le c_1 \frac{\phi(r^{-2})^{1/2}}{\phi(\delta_{\overline{B}(0,r)^c}(y)^{-2})^{1/2}} g(r)\le c_2 {\phi(r^{-2})^{-1/2}}{\phi(\delta_{\overline{B}(0,r)^c}(y)^{-2})^{-1/2}} r^{-d} .
\end{equation}
Since $G_{\overline{B}(0,r)^c}(x,y) \le g(|x-y|)$, by \eqref{e:G} and \eqref{e:newge}, we only need to prove the first inequality in \eqref{e:newge2}
for $(x, y)$ satisfying $y \in A(0,r,(p+7)r/8)$ and $|x-z|\le 2|x-y|$.
In this case, we have
\begin{equation} \label{e:newge3}
|x-y| \ge \frac12(p-1)r >4\delta_{\overline{B}(0,r)^c}(y).
\end{equation}
Let $y_1:=   (8^{-1}(p-1)+1)  r|y|^{-1}  y$.
Then
$$
|y_1-y| \le \delta_{\overline{B}(0,r)^c}(y) \vee \delta_{\overline{B}(0,r)^c}(y_1) \le \frac18 (p-1)r \le  \frac14 |x-y|.$$
Thus by \eqref{e:newge3}
\begin{equation} \label{e:newge4}
|x-y_1| \ge |x-y| -|y_1-y| \ge \frac{3}{4}|x-y| \ge  \frac{3}{8}(p-1)r.
\end{equation}
 Because of  \eqref{e:newge3}, we
can apply Theorem \ref{L:2} and then
use  \eqref{e:newge4} to get that for
$(x, y)$ satisfying $y \in A(0,r,(p+7)r/8)$ and $|x-z|\le 2|x-y|$,
\begin{eqnarray*}
   && G_{\overline{B}(0,r)^c}(x,y)
    \le
     c_3\frac{\phi(\delta_{\overline{B}(0,r)^c}(y_1)^{-2})^{1/2}}{\phi(\delta_{\overline{B}(0,r)^c}(y)^{-2})^{1/2}}\,  G_{\overline{B}(0,r)^c}(x,y_1)\le
     c_3\frac{\phi((64) (p-1)^{-2} r^{-2})^{1/2}}{\phi(\delta_{\overline{B}(0,r)^c}(y)^{-2})^{1/2}}\, g(|x-y_1|)
\\
&&     \le
     c_4\frac{\phi(r^{-2})^{1/2}}{\phi(\delta_{\overline{B}(0,r)^c}(y)^{-2})^{1/2}}\, g(\frac{3}{8}(p-1)r)
      \le
     c_5\frac{\phi(r^{-2})^{1/2}}{\phi(\delta_{\overline{B}(0,r)^c}(y)^{-2})^{1/2}}\, g(r)
     \,
\end{eqnarray*}
with  constants $c_i=c_i(\phi,p,q)>0$, $i=3,4,5$.
In the last inequality we have used \eqref{e:doubling-conditionG}.
Therefore using \eqref{e:G}  we have proved \eqref{e:newge2}.

Applying \eqref{e:newge2} to $I_1$ and using the fact that $\delta_{\overline{B}(0,r)^c}(y) \le |y-z|$,
 we get
$$
    I_1 \le  c_2  r^{-d}  {\phi(r^{-2})^{-1/2}} \int_{A(0,r,2qr)\cap\{|x-z|> 2|x-y|\}}
   \phi(|y-z|^{-2})^{-1/2}j(|y-z|)\, dy\, .
$$
Since $B(0,r)^c\subset \overline{B}(z,r-|z|)^c$,
by \eqref{e:J}, the integral above is less than or equal to
\begin{eqnarray*}
       \int_{\overline{B}(z,r-|z|)^c} \phi(|y-z|^{-2})^{-1/2} j(|y-z|)\,\, dy
    \,\le \,  c_6 \int_{r-|z|}^{\infty} \frac{\phi(t^{-2})^{1/2}}{t}\, dt
  \, \le \, c_7\phi((r-|z|)^{-2})^{1/2}\, ,
\end{eqnarray*}
where in the last inequality we used \eqref{e:ie-3}.
Hence, $I_1\le c_{8} \phi(r^{-2})^{-1/2} \phi((r-|z|)^{-2})^{1/2} r^{-d}$.

To estimate $I_2$ we first
note that  if $2|x-y|\le |x-z|$, then $|y-z|\ge |x-z|-|y-x|\ge \frac12 |x-z|$,
hence, by \eqref{c:doubling-condition},  $j(|y-z|)\le c_{9} j(|x-z|)$ where $c_{9}=c_{9}(\phi)>0$. Thus,
\begin{align*}
      & I_2    \,\le \, c_{10}j(|x-z|) \int_{B(x,\frac{|x-z|}{2})}
      g(|x-y|)\, dy
    \,\le \, c_{11}j(|x-z|)\int_0^{\frac{|x-z|}{2}} t^{-1}\phi(t^{-2})^{-1}\, dt \\
    &\le \, c_{12}j(|x-z|)\phi(|x-z|^{-2})^{-1}
   \,\le\, c_{13}|x-z|^{-d} \le c_{14} r^{-d}\, .
\end{align*}
In the penultimate inequality,  we used \eqref{e:ie-3}.

Finally, we deal with $I_3$. For $|y| \ge 2qr$ and $|z|<r$
we have that $|y-z|\ge |y|-|z|>|y|-r\ge (1-1/(2q))|y|$, hence by \eqref{e:doubling-condition}, we get $j(|y-z|)\le c_{15} j(|y|)$. Also, for $|x|<qr$ and $|y|\ge 2qr$, we have that $|y-x|\ge |y|-|x|\ge |y|/2$, hence $
g(|x-y|)\le c_{16}
g(|y|)$ by \eqref{e:doubling-conditionG}.
Therefore, by Theorem \ref{t:J-G},
\begin{align*}
  I_3\,
\le \, c_{17} \int_{\overline{B}(0,2qr)^c}  \frac{1}{|y|^d \phi(|y|^{-2})} \, \frac{\phi(|y|^{-2})}{|y|^d}\, dy
    \,\le\,  c_{18} \int_{2qr}^{\infty} t^{-d-1}\, dt
   \,=\,c_{19} r^{-d}\, .
\end{align*}
This concludes the proof of the lemma. \qed

It is easy to see that, by the strong Markov property,
for all Greenian open sets $U$ and $D$ with $U \subset D$,
$
G_D(x,y)= G_U(x,y) + \E_x\left[   G_D(X_{\tau_U}, y)\right]
$ for every $(x,y) \in \R^d \times \R^d.
$
Thus,
for all Greenian open sets $U$ and $D$ with $U \subset D$,
\beq\label{e:KDU}
K_D(x,z)= K_U(x,z) + \E_x\left[   K_D(X_{\tau_U}, z)\right], \qquad (x,z) \in
U \times \overline{D}^c.
\eeq

Since $X$ is a purely discontinuous rotationally invariant L\'evy process,
it follows from \cite[Proposition 4.1]{Mi}
(see also \cite[Theorem 1]{Sz1}) that
if $V$ is a Lipschitz open set and $U \subset V$,
\beq\label{e:Lip}
\P_x(X_{\tau_U} \in \partial V)=0 \qquad \text{ and } \qquad \P_x(X_{\tau_U} \in dz) = K_U(x,z)dz \quad\text{on }V^c.
\eeq

\begin{lemma}\label{l:sup-estimate}
Let $1<p<q<\infty $.
There exists $c=c(\phi,p,q)>1$ such that for all
$r\ge 1/2$ and all open sets $U\subset \overline{B}(0,r)^c$ it holds that
\begin{equation}\label{e:sup-estimate}
K_U(x,y)\le c r^{-d} \left(\int_{U\cap B(0, \frac{1+p}{2}r)} K_U(z,y)\, dz +1\right) ,\quad \textrm{for all }
x\in A(0,pr,qr)\cap U, \, y\in B(0,r)\, .
\end{equation}
\end{lemma}
\pf Let $q_1=\frac{3+p}{4} $, $q_2=\frac{1+p}{2}$ (so that $1<q_1<q_2<p$) and $s\in [q_1 r, q_2 r]$.
Then, by \eqref{e:KDU} and \eqref{e:Lip},
for $ x\in A(0,pr,qr)\cap U$ and $ y\in B(0,r)$ it holds that
\begin{eqnarray*}
    K_U(x,y) &=& \E_x\left[K_U(X_{\tau_{U\cap \overline{B}(0,s)^c}}, y)\right] + K_{U\cap \overline{B}(0,s)^c}(x,y)\\
    & = & \int_{U\cap B(0,s)}K_U(z,y) K_{U\cap \overline{B}(0,s)^c}(x,z)\, dz + K_{U\cap \overline{B}(0,s)^c}(x,y)\\
    &\le & \int_{U\cap B(0,q_2r)}\ind_{\{|z|<s\}} K_U(z,y) K_{\overline{B}(0,s)^c}(x,z)\, dz  + K_{\overline{B}(0,s)^c}(x,y).
\end{eqnarray*}
Hence by Fubini's theorem,
\begin{align*}
  &  (q_2-q_1)r K_U(x,y) \,=\,  \int_{q_1 r}^{q_2 r} K_U(x,y)\, ds\\
    &\le   \int_{U\cap B(0,q_2 r)}\left(\int_{|z|}^{q_2 r} K_{\overline{B}(0,s)^c}(x,z)\, ds\right)\, K_U(z,y)\, dz+ \int_{q_1 r}^{q_2 r} K_{\overline{B}(0,s)^c}(x,y)\, ds
    \,=:\, I_1+I_2\, .
\end{align*}

For $s\in [q_1r, q_2r]$ and $z\in U\cap B(0, s)$, we have $r<|z|<s\le q_2 r = \frac{1+p}{2} r <pr<|x| \le qr \le (q/q_1)s$, so
 it follows from
Lemmas \ref{l:phi-property} and \ref{l:poisson-kernel-estimate} that
\begin{eqnarray*}
    K_{\overline{B}(0,s)^c}(x,z) &\le & c_1 \Big(
    s^{-d}\big(\phi(s^{-2})^{-1/2} \phi((s-|z|)^{-2})^{1/2}+1\big)+s^{-d}\Big)\\
    & \le & c_2 \Big(r^{-d}\big(\phi((q_2 r)^{-2})^{-1/2}\phi((s-|z|)^{-2})^{1/2}+1\big)+r^{-d}\Big)\\
    & \le & c_3 r^{-d}\Big(\phi(r^{-2})^{-1/2} \phi((s-|z|)^{-2})^{1/2}+1\Big )\, ,
\end{eqnarray*}
where $c_2=c_2(\phi, p, q)$ and $c_3=c_3(\phi, p, q)$. Hence,
\begin{eqnarray*}
    \int_{|z|}^{q_2 r} K_{\overline{B}(0,s)^c}(x,z)\, ds
    &\le &c_3 r^{-d}\Big(\phi(r^{-2})^{-1/2} \int_0^{q_2 r-|z|}\phi(t^{-2})^{1/2}\, dt + q_2 r\Big)\\
    &\le & c_4 r^{-d} \Big(\phi(r^{-2})^{-1/2} (q_2 r- |z|) \phi((q_2 r-|z|)^{-2})^{1/2} +r\Big)\, ,
\end{eqnarray*}
where the last inequality follows from \eqref{e:ie-1}.
Since $t\mapsto t\phi(t^{-2})^{1/2}$ is increasing by \eqref{e:Berall}
and $q_2r-|z|\le q_2 r$, we have that $(q_2 r- |z|) \phi((q_2 r-|z|)^{-2})^{1/2}\le q_2r \phi((q_2 r)^{-2})^{1/2}\le c_{5} r \phi(r^{-2})^{1/2}$. Therefore,
$$
    \int_{|z|}^{q_2 r} K_{\overline{B}(0,s)^c}(x,z)\, ds \le c_6 r^{-d}\Big(\phi(r^{-2})^{-1/2} r \phi(r^{-2})^{1/2} +r\Big)\le c_7 r^{-d+1}\, .
$$
Further, for $s\in [q_1r, q_2r]$ and $y\in B(0, r)$ we have $|y|<r<q_1r <s$, and so $s-|y|>(q_1-1)r$,
and we get similarly as above (but easier) that
$$
    I_2=\int_{q_1 r}^{q_2 r} K_{\overline{B}(0,s)^c}(x,y)\, ds \le c_7 r^{-d+1}\, .
$$
Finally,
\begin{eqnarray*}
    K_U(x,y)= \frac{1}{(q_2-q_1)r}(I_1+I_2) \le  \frac{p-1}{4} c_8 r^{-d}\left(\int_{U\cap B(0,q_2 r)}K_U(z,y)\, dz + 1\right)\, ,
\end{eqnarray*}
proving the lemma. \qed

\begin{lemma}\label{l:key-lemma}
For every $a>1$ there exists $c=c(\phi, a)>1$ such that for all $r\ge 1$ and all open sets $U\subset \overline{B}(0,r)^c$ it holds that
\begin{align}\label{e:key-lemma}
\frac{1}{c}\, K_U(x,0)\left(\int_{U\cap B(0,
ar
)}K_U(y,z)\, dy +1\right)\le K_U(x,z)
   \le  c K_U(x,0)\left(\int_{U\cap B(0,
ar
)}K_U(y,z)\, dy +1\right)
\end{align}
for all  $x\in U\cap \overline{B}(0,
ar
)^c$ and $z\in B(0,r)$.
\end{lemma}

\pf
Fix two constant $b_2$ and $b_3$ such that $1>b_2>b_3>1/a$. Let
$U_1:=\overline{B}(0, ar)^c\cap U$, $U_2:=\overline{B}(0, b_2ar)^c\cap U$ and
$U_3:=\overline{B}(0, b_3ar)^c\cap U$.
Then by \eqref{e:KDU},
for $x\in U_1$ and $z\in B(0,r)$,
\begin{eqnarray*}
K_U(x,z)&=&\E_x\left[K_U(X_{\tau_{U_2}},z)\right] +K_{U_2}(x,z)\\
    &=&\int_{U_3\setminus U_2}K_U(y,z)\P_x(X_{\tau_{U_2}}\in dy)+\int_{U\setminus U_3}K_U(y,z)K_{U_2}(x,y)\, dy +K_{U_2}(x,z)\\
    &=:& I_1+I_2+I_3\, .
\end{eqnarray*}

We first estimate $I_3=K_{U_2}(x,z)=\int_{U_2}G_{U_2}(x,y)j(|y-z|)\, dy$.
For $y\in U_2$ we have that $|y|> b_2ar$, hence
$$
(1-\frac1{b_2a})|y|\le |y|-|z|\le |y-z|\le |y|+|z|\le (1+\frac1{b_2a})|y|\, .
$$
Hence, by \eqref{e:doubling-condition}, there exists
$c_1=c_1(\phi)>0$  such that $c_1^{-1}j(|y|)\le j(|y-z|)\le c_1 j(|y|)$. Therefore,
\begin{eqnarray}
    c_1^{-1} K_{U_2}(x,0)&=&c_1^{-1}\int_{U_2}G_{U_2}(x,y)j(|y|)\, dy \le \int_{U_2}G_{U_2}(x,y)j(|y-z|)\, dy \label{e:I3-lower}\\
    &= & K_{U_2}(x,z)\le c_1 K_{U_2}(x,0)\le c_1 K_U(x,0)\, . \label{e:I3-upper}
\end{eqnarray}

In order to estimate $I_2=\int_{U\setminus U_3}K_U(y,z)K_{U_2}(x,y)\, dy$, we proceed similarly
by estimating $K_{U_2}(x,y)= \int_{U_2}G_{U_2}(x,w)j(|w-y|)\, dw$ for
$x\in U\cap \overline{B}(0, ar)^c$ and $y\in U\setminus U_3$.
Note that since $y\notin U_3$ it holds that $|y|< b_3ar$.
For $w\in U_2$ it holds that $|w|>b_2ar$.
Hence, similarly as above we get that $(1-\frac{b_3}{b_2}) |w|\le |w-y|\le (1+\frac{b_3}{b_2}) |w|$.
Thus there exists $c_2=c_2(\phi)>0$  such that $c_2^{-1}j(|w|)\le j(|w-y|)\le c_2 j(|w|)$.
In the same way as above, this implies that
$c_2^{-1}K_{U_2}(x,0)\le K_{U_2}(x,y)\le c_2 K_{U_2}(x,0)\le c_2 K_U(x,0)$. Therefore
\begin{align}
    \lefteqn{c_2^{-1} K_{U_2}(x,0)\int_{U\setminus U_3} K_U(y,z)\, dy \le I_2}\label{e:I2-lower}\\
    &\le \,c_2 K_{U_2}(x,0)\int_{U\setminus U_3} K_U(y,z)\, dy \le c_2 K_{U}(x,0)\int_{U\cap B(0,
ar
)} K_U(y,z)\, dy\, . \label{e:I2-upper}
\end{align}

In the case of $I_1$ we only need an upper estimate. It holds that
\begin{eqnarray}
    I_1=\int_{U_3\setminus U_2}K_U(y,z)\P_x(X_{\tau_{U_2}}\in dy)
    \le  \left(\sup_{y\in U_3\setminus U_2}K_U(y,z)\right)
    \P_x\left(X_{\tau_{U_2}}\in \overline{B}(0, b_2ar)\right)\, .
    \label{e:I1}
\end{eqnarray}
By Lemma \ref{l:sup-estimate} (with $p=b_3a$ and $q=b_2a$),
there is a constant $c_3=c_3(\phi, a)>0$ such that
$$    \left(\sup_{y\in U_3\setminus U_2}K_U(y,z)\right) \le c_3 r^{-d} \left(\int_{U\setminus U_3}K_U(y,z)\, dy +1\right) \, .
$$
By Lemma \ref{l-exit} used with $D=U_2$, $b_2ar$ instead of $r$ and
$\frac1{b_2}$ instead of $a$, there is a constant
$c_4=c_4(\phi, a)>0$ such that
$$
\P_x\left(X_{\tau_{U_2}}\in \overline{B}(0, b_2ar)\right)
\le c_4 r^d K_{U_2}(x,0)\, .
$$
By applying the last two estimates to \eqref{e:I1} we get
\begin{equation}\label{e:I1-b}
    I_1 \le c_5 K_{U_2}(x,0) \left( \int_{U\setminus U_3}K_U(y,z) \, dy +1\right)\, .
\end{equation}
Putting together \eqref{e:I3-upper}, \eqref{e:I2-upper} and \eqref{e:I1-b}, we see that
\begin{eqnarray}
    K_U(x,z)& \le & c_6 K_{U_2}(x,0)\left(\int_{U\setminus U_3}K_U(y,z)\, dy +1 \right) \label{e:all-upper}\\
    & \le & c_6 K_U(x,0)\left(\int_{U\cap B(0,
ar
)}K_U(y,z)\, dy +1 \right)\, .\nonumber
\end{eqnarray}
Thus, the upper bound in \eqref
{e:key-lemma} holds true.

In order to prove the lower bound, we may neglect $I_1$. First we note that for $z\in B(0,r)$,
\begin{eqnarray}
    \int_{U\cap B(0,
ar
)}K_U(y,z)\, dy &=& \int_{U\setminus U_3}K_U(y,z)\, dy +\int_{U_3\setminus U_1} K_U(y,z)\, dy \nonumber \\
    & \le & \int_{U\setminus U_3}K_U(y,z)\, dy +\left(\sup_{y\in U_3\setminus U_1} K_U(y,z)\right) \big| U_3\setminus U_1\big| \nonumber \\
    & \le & \int_{U\setminus U_3}K_U(y,z)\, dy +c_3 r^{-d} \left(\int_{U\setminus U_3}K_U(y,z)\, dy +1\right) c_7 r^d \nonumber \\
    & \le & c_8 \left(\int_{U\setminus U_3}K_U(y,z)\, dy +1\right)\, . \label{e:comparison-of-integrals}
\end{eqnarray}
Here we used Lemma \ref{l:sup-estimate} in the second inequality.
Next, by using the already proved upper bound \eqref{e:all-upper} with $z=0$, we see that
\begin{eqnarray}
    K_{U}(x,0) & \le  & c_6 K_{U_2}(x,0)\left(\int_{U\setminus U_3}K_U(y,0)\, dy +1 \right) \nonumber\\
    & \le  & c_6 K_{U_2}(x,0)\left(\int_{
    A(0, r, 3r/2) }K_{\overline{B}(0,\frac{r}{2})^c}(y,0)\, dy +1 \right)\nonumber \\
      & \le & c_6 K_{U_2}(x,0)\left(\int_{
       A(0, r, 3r/2)} c_9 r^{-d}\, dr+1\right)\, \le \, c_{9}\, K_{U_2}(x,0)\, . \label{e:comparison-U2-U}
\end{eqnarray}
Here we have used Lemma \ref{l:poisson-kernel-estimate} in the
third  inequality.
The lower estimate now follows from \eqref{e:I3-lower}, \eqref{e:I2-lower}, \eqref{e:comparison-of-integrals} and \eqref{e:comparison-U2-U}.
\qed

\pff {\bf of Theorem \ref{bhp-inf}.} Let $x\in U\cap \overline{B}(0,
ar)^c$.
Then, by \eqref{e:Lip},
\begin{eqnarray*}
    u(x)&=& \int_{B(0,r)}K_U(x,z)u(z)\, dz\\
    &\asymp &  K_U(x,0)\int_{B(0,r)}\left(\int_{U\cap B(0, ar
)}K_U(y,z)\, dy +1\right)u(z)\, dz\\
    &=& K_U(x,0)\left(\int_{B(0,r)}u(z)\, dz +\int_{U\cap B(0, ar
)}\left(\int_{B(0,r)}K_U(y,z)u(z)\, dz\right)\, dy\right)\\
    &=& K_U(x,0)\left(\int_{B(0,r)}u(z)\, dz +\int_{U\cap B(0, ar
)}u(y)\, dy\right)\\
    &\asymp& K_U(x,0)\int_{B(0, ar)}u(z)\, dz\, ,
\end{eqnarray*}
where in the second line we used  Lemma \ref{l:key-lemma} and in the first, fourth and last line we used the fact that $u$ vanishes a.e. on $\overline{B}(0,r)^c\setminus U$.
\qed

\pff {\bf of Corollary \ref{c:bhp-inf}.}
It follows from \eqref{e:bhp-inf} that for $x,y\in U\cap \overline{B}(0, ar )^c$,
$$
    \frac{u(x)}{v(x)} \le \frac{{C_1} K_U(x,0)\int_{B(0, ar)}u(z)\, dz}{C_1^{-1}K_U(x,0)\int_{B(0,ar)}v(z)\, dz} =
C_1^2 \, \frac{\int_{B(0,ar)}u(z)\, dz}{\int_{B(0,ar)}v(z)\, dz}\, .
$$
Similarly,
$$
\frac{u(y)}{v(y)} \ge \frac{C_1^{-1} K_U(y,0)\int_{B(0,ar)}u(z)\, dz}{C_1K_U(y,0)\int_{B(0,ar)}v(z)\, dz} =
C_1^{-2}\, \frac{\int_{B(0,ar)}u(z)\, dz}{\int_{B(0,ar)}v(z)\, dz}\, .
$$The last two displays show that \eqref{e:bhp-inf-cor} is true for $x,y\in U\cap \overline{B}(0,
ar)^c$ with $C_2=C_1^4$.
\qed

\begin{corollary}\label{c:hi}
For every $a>1$,  there exists $c=c(d, \phi, a)>1$  such that
\begin{itemize}
\item[(i)]
    for every $r\ge 1$, every open set $U\subset \overline{B}(0,r)^c$ and every nonnegative function $u$ on
    $\R^d$ which is regular harmonic in $U$ and vanishes  a.e on $\overline{B}(0,r)^c\setminus U$,
    $$
    \frac{u(x)}{K_U(x,0)}\le
c\frac{u(y)}{K_U(y,0)}\, ,\qquad \textrm{for all }x,y\in U\cap \overline{B}(0,ar)^c
\,;
    $$
\item[(ii)]
for every $r\ge 1$ and every nonnegative function $u$ on $\R^d$ which is
regular harmonic in $\overline{B}(0,r)^c$,
    $$
    \frac{u(x)}{K_{\overline{B}(0,r)^c}(x,0)}\le
c
 \frac{u(y)}{K_{\overline{B}(0,r)^c}(y,0)}\, ,\qquad \textrm{for all }x,y\in \overline{B}(0,
ar
)^c\, .
    $$
\end{itemize}
\end{corollary}
\pf The first claim is a direct consequence of Theorem \ref{bhp-inf} with $
c
=C_1^2$, while the second follows from the first and the fact that the zero boundary condition is vacuous. \qed

\begin{lemma}\label{l:exit-probability}
For every $a>1$, there exists $c=c(d,\phi, a)>1$ such that for all $r>0$
$$
c^{-1}r^d \phi(r^{-2}) G(x,0) \le \P_x(\tau_{\overline{B}(0,r)^c}<\infty) \le
c r^d \phi(r^{-2}) G(x,0)\, ,\quad \textrm{for all } x\in \overline{B}(0,
ar
)^c\, .
$$
\end{lemma}
\pf Let $x\in \overline{B}(0,
ar
)^c$. By the strong Markov property,
\begin{equation}\label{e:smp}
    \int_{\overline{B}(0,r)} G(x,y)\, dy=\E_x\left[\int_{\overline{B}(0,r)}G(X_{\tau_{\overline{B}(0,r)^c}},y)\, dy \, , \tau_{\overline{B}(0,r)^c}<\infty\right].
\end{equation}
Since $r\mapsto
g(r)$ is decreasing, \cite[Lemma 5.53]{BBKRSV} shows that there exists a constant $c=c(d)$ such that for every $r>0$ and all $z\in \overline{B}(0,r)$ we have
$$
    c \int_{
    \overline{B}(0,r)}
g(|y|)\, dy \le \int_{
    \overline{B}(0,r)} G(z,y)\, dy \le \int_{
    \overline{B}(0,r)}
g(|y|)\, dy \, .
$$
Then it follows from \eqref{e:smp} that for $x\in \overline{B}(0,
ar)^c $,
$$
    \int_{\overline{B}(0,r)} G(x,y)\, dy \asymp \left(\int_{B(0,r)}
g(|y|)\, dy \right) \P_x(\tau_{\overline{B}(0,r)^c}<\infty)\, ,
$$
with a constant depending on $d$ only. By the uniform  Harnack inequality (Theorem \ref{t:uhi}), there exists $c_2>1$ such that
$
    c_2^{-1} G(x,0)\le G(x,y) \le c_2 G(x,0)$
  for every $ x\in \overline{B}(0,
ar
)^c$ and  $ y\in B(0,r)$.
Hence
\begin{align}
    r^d G(x,0)\asymp \left(\int_{B(0,r)}
g(|y|)\, dy \right) \P_x(\tau_{\overline{B}(0,r)^c}<\infty)\, , \quad x\in \overline{B}(0,ar)^c,
\label{e:new34}
\end{align}
with a constant depending on
$d$ and $a$. It follows from
\eqref{e:G} and  \eqref{e:ie-3}  that
\begin{align}
    \int_{B(0,r)}
g(|y|)\, dy \asymp \frac{1}{\phi(r^{-2})}\, ,
\label{e:new35}
\end{align}
with a constant depending on
$d$ and $a$.
Combining \eqref{e:new34}-\eqref{e:new35} we have proved the lemma
\qed

\begin{corollary}\label{c:poisson-inf}
For every $a>1$, there exists $c=c(\phi,d, a)>1$ such that for all $r\ge 1$ it holds that
\begin{equation}\label{e:G-K}
c^{-1} \phi(r^{-2}) G(x,0)\le K_{\overline{B}(0,r)^c}(x,0)\le
c\phi(r^{-2}) G(x,0)\, ,\quad x\in \overline{B}(0, ar
)^c\, ,
\end{equation}
and consequently
$$
    \lim_{|x|\to \infty}K_{\overline{B}(0,r)^c}(x,0)=0\, .
$$
\end{corollary}
\pf Note that $\P_x(\tau_{\overline{B}(0,r)^c}<\infty)=\int_{B(0,r)}K_{\overline{B}(0,r)^c}(x,z)\, dz$. Further, for $z\in B(0,r)$ and $y\in \overline{B}(0,r)^c$ we have that $|y-z|\le 2|y|$ and hence $j(|y-z|)\ge c_1 j(|y|)$ by \eqref{e:doubling-condition}. Therefore,
\begin{eqnarray*}
K_{\overline{B}(0,r)^c}(x,z)
\,\ge\, c_1\, \int_{\overline{B}(0,r)^c}G_{\overline{B}(0,r)^c}(x,y)j(|y|)\, dy \,=\,c_1\, K_{\overline{B}(0,r)^c}(x,0)\, .
\end{eqnarray*}
Using \eqref{e:Lip} it follows that
$\P_x(\tau_{\overline{B}(0,r)^c}<\infty)\ge c_1 \int_{B(0,r)}K_{\overline{B}(0,r)^c}(x,0)\, dz =c_2 r^d K_{\overline{B}(0,r)^c}(x,0)$. On the other hand, from Lemma \ref{l-exit} with $D=\overline{B}(0,r)^c$, we see that $\P_x(\tau_{\overline{B}(0,r)^c}<\infty)\le c_2 r^d K_{\overline{B}(0,r)^c}(x,0)$. Thus
$$
    \P_x(\tau_{\overline{B}(0,r)^c}<\infty)\asymp r^d K_{\overline{B}(0,r)^c}(x,0)\, , \quad x\in \overline{B}(0,
ar
)^c\, .
$$
Comparing with the result in Lemma \ref{l:exit-probability} gives \eqref{e:G-K}. The last statement
follows from $\lim_{|x|\to \infty}G(x, 0)=0$.
\qed

\begin{corollary}\label{c:vanish-at-infinity}
Let $r\ge 1$ and $U\subset \overline{B}(0,r)^c$. If $u$ is a non-negative function on $\R^d$ which is regular harmonic in $U$ and
vanishes a.e.~on  $\overline{B}(0,r)^c\setminus U$, then
$$
    \lim_{|x|\to \infty}u(x)=0\, .
$$
\end{corollary}
\pf Note that $K_U(x,0)\le K_{\overline{B}(0,r)^c}(x,0)$. It follows from Corollary \ref{c:poisson-inf} that  $$\lim_{|x|\to \infty}K_{\overline{B}(0,r)^c}(x,0)=0.$$ Then the claim follows from Theorem \ref{bhp-inf}. \qed

\begin{remark}
{\rm
(i)
Corollary \ref{c:vanish-at-infinity} is not true if \emph{regular harmonic} is replaced by \emph{harmonic}. Indeed,
let $V$ denote the renewal function of the ladder height process of the one-dimensional subordinate Brownian motion $W^d(S_t)$.
Then the function $w(x)=w(\widetilde{x}, x_d):=V((x_d)^+)$ is harmonic in the
upper half-space $\bH\subset B((\widetilde{0},-1), 1)^c$ (see \cite{KSV3}),
vanishes on $B((\widetilde{0},-1), 1)^c \setminus \bH$, but clearly $\lim_{x_d\to \infty}w(x)=\infty$.

\noindent
(ii)
When $d=1\le \alpha$, Corollary \ref{c:vanish-at-infinity} is not true even for the symmetric $\alpha$-stable
process
because the Green function of the complement
of any bounded interval does not vanish at infinity, which can be seen using the Kelvin
transform.}
\end{remark}

%%%%%%%%%%%%%%%%%%%%%%%%%%%%%%%%%%%%%%%%%%%%%%%%%%%%  The infinite part of the Martin boundary  %%%%%%%%%%%%%%%%%%%%%%%%%%%%%%%%%%%%%

\section{The infinite part of the Martin boundary}

In this section we will consider a large class of
unbounded open sets $D$ and identify the infinite part of the
Martin boundary of $D$ without assuming that the finite part of the Martin boundary of $D$ coincides with the
Euclidean boundary.

We first recall the definition of $\kappa$-fatness  at infinity from the introduction: Let $\kappa \in (0,1/2]$. An open set $D$ in $\R^d$ is $\kappa$-fat at infinity if there exists $R>0$ such that for every $r\in [R,\infty)$ there exists $A_r \in \R^d$ such that $B(A_r, \kappa r)\subset D\cap \overline{B}(0,r)^c$ and $|A_r|< \kappa^{-1} r$.

The origin does not play any special role in this definition:
Suppose that $D$ is $\kappa$-fat at infinity with characteristics $(R,\kappa)$. For every $Q\in \R^d$, define $R_Q:=R \vee |Q|$.
For all $r\ge R_Q$, with $\wh A_r:=A_{2r}$ and $\wh \kappa:=\kappa/3$, we have
$$
B(\wh A_r, \wh \kappa r)\subset B(A_{2r}, 2\kappa r)\subset D\cap \overline{B}(0,2r)^c \subset D\cap \overline{B}(0,r+|Q|)^c\subset D\cap \overline{B}(Q,r)^c
$$
and
$
|\wh A_r-Q|\le|A_{2r}|+|Q| \le (\kappa/2)^{-1}r+r   < (\wh \kappa)^{-1}r.
$

\begin{remark}\label{r:kappa-fat-infty}
{\rm
(i) Note that it follows from the definition that any
open set which is $\kappa$-fat at infinity is necessarily unbounded.

\noindent
(ii) Since $B(A_r, \kappa r)\subset \overline{B}(0,r)^c$ we have that $|A_r|-\kappa r >r$
implying $(\kappa+1)r <|A_r| <\kappa^{-1} r$.

\noindent
(iii) We further note that $B(A_r, (\kappa/2) r)\cap B(A_{(\kappa/2)^{-1}r},  r)=\emptyset$.
Indeed, for any point $x$ in the intersection we would have that  $|x|\le |A_r|+(\kappa/2) r
< (\kappa^{-1}+\kappa/2) r$, and at the same time $|x|\ge |A_{(\kappa/2)^{-1}r}|- r>
(\kappa+1)(\kappa/2)^{-1}r-\kappa r=(2\kappa^{-1}+2-\kappa) r$. But this is impossible.
}
\end{remark}

In this section we first identify the infinite part of the Martin boundary of an open set
$D\subset \R^d$ which is $\kappa$-fat at infinity with characteristics $(R,\kappa)$.
Without loss of generality, we assume that $R>1$.

Recall that,  throughout the paper we assume that {\bf (H1)} and {\bf (H2)}
are true and $d >2(\delta_2 \vee \delta_4)$.

\begin{lemma}\label{l:lemma3.1}
Let $D\subset \R^d$ be an open set which is $\kappa$-fat at infinity with
characteristics $(R,\kappa)$. There exist
$c=c(d,\phi, \kappa)>0$ and $\gamma=\gamma(d,\phi, \kappa)\in (0,d)$ such that for every $r\ge R$
and any non-negative function $h$ in $\R^d$ which is harmonic in $D\cap \overline{B}(0,r)^c$ it holds that
\begin{equation}\label{e:lemma3.1}
    h(A_r)\le c (\kappa/2)^{-(d-\gamma)k}h(A_{(\kappa/2)^{-k} r})\, ,\quad k=0,1,2,\dots \, .
\end{equation}
\end{lemma}
\pf Fix $r\ge R$. For $n=0,1,2,\dots $,
let $\eta_n=(\kappa/2)^{-n} r$,  $A_n=A_{\eta_n}$ and $B_n=B(A_n,  \eta_{n-1})$ (where $\eta_{-1}= (\kappa/2) r$). Note that the balls $B_n$ are pairwise disjoint (cf.~Remark \ref{r:kappa-fat-infty} (iii)). By harmonicity of $h$, for every $n=0,1,2 \dots $,
\begin{eqnarray*}
    h(A_n)
    = \E_{A_n}\left[h(X_{\tau_{B_n}})\right]
    \ge \sum_{l=0}^{n-1} \E_{A_n}\left[h(X_{\tau_{B_n}}):\, X_{\tau_{B_n}}\in B_l\right]
   = \sum_{l=0}^{n-1} \int_{B_l} K_{B_n}(A_n, z)h(z)\, dz\, .
\end{eqnarray*}
By the uniform Harnack inequality, Theorem \ref{t:uhi}, there exists
$c_1=c_1(d, \kappa, \phi)>0$
 such that for every $l=0,1,2,\dots$,
$
    h(z)\ge c_1 h(A_l)$ for all  $ z\in B_l$ .
Hence
$$
      \int_{B_l}K_{B_n}(A_n, z)h(z)\, dz \ge c_1 h(A_l) \int_{B_l}K_{B_n}(A_n, z)\, dz\, , \quad
      0\le l\le n-1\, .
$$
By \eqref{e:pke-lower} we have
$$
  \int_{B_l}  K_{B_n}(A_n,z)dz\ge c_2\phi(\eta_n^{-2})^{-1}\int_{B_l} j(|2(A_n-z)|)dz\, , \quad  0\le l\le n-1\, .
$$
For $z\in B_l$, $l=0, 1, \cdots, n-1$, it holds that $|z|\le \kappa^{-1}(\kappa/2)^{-l}r
+(\kappa/2)^{-(l-1)}r
=(\kappa/2)^{-l}r(\kappa^{-1}+\kappa/2)$.
Since $ |A_n|\le \kappa^{-1}  \eta_n  $, we have  that $|A_n-z|\le  |A_n|+|z|\le 2\kappa^{-1}\eta_n$.
Together with Theorem \ref{t:J-G} and Lemma \ref{l:phi-property}, this implies that
$
    j(|2(A_n-z)|)\ge c_3 j(|\eta_n|)$ for every  $z\in B_l$ and  $0\le l\le n-1$.
Therefore,
$$
    \int_{B_l}K_{B_n}(A_n, z)\, dz \ge c_4 j(|\eta_n|)\phi(\eta_n^{-2})^{-1} |{B_l}| \ge c_5\eta_n^{-d}\eta_l^{d} =c_5 \frac{\eta_l^d}{\eta_n^d}
    \, , \quad 0\le l\le n-1\, .
$$
Hence,
$$
    \eta_n^d h(A_n)\ge c_5 \sum_{l=0}^{n-1}\eta_l^d h(A_l)\, ,\quad \textrm{for all }n=1,2,\dots .
$$
Let $a_n:=\eta_n^d h(A_n)$ so that $a_n\ge c_5 \sum_{l=0}^{n-1} a_l$.
Using the identity $1+ c_5+ c_5 \sum_{l=1}^{n-2} (1+c_5)^{l}= (1+c_5)^{n-1}$ for $n \ge3$,
by induction it follows that $a_n\ge c_5 (1+c_5)^{n-1} a_0$.
 Let $\gamma:=\log\left(1+
c_5\right)/\log(2/\kappa)$ so that $\left(1+
c_5\right)^n = (2/\kappa)^{\gamma n}$. Note that $c_5$ can be chosen arbitrarily
close to  zero (but positive), so that $\gamma  <d$. Thus, $a_0\le
(1+c_5)c_5^{-1} (\kappa/2)^{\gamma n} a_n$, or
$
    \eta_0^d h(A_0)\le
(1+c_5)c_5^{-1} (\kappa/2)^{\gamma n} \eta_n^d h(A_n)\, .
$
Hence,
$
    h(A_r)\le
(1+c_5)c_5^{-1} (\kappa/2)^{\gamma n} (\kappa/2)^{-d n} h(
A_{(\kappa/2)^{-n} r}
)\, .
$
\qed

\begin{lemma}\label{l:lemma3.2}
Let $D\subset \R^d$ be an open set which is $\kappa$-fat at infinity with
characteristics $(R,\kappa)$. There exists $c=c(d,\phi, \kappa)>0$ such that for every
$r \ge R$ and every non-negative function $h$ on $\R^d$
which is regular harmonic in $D\cap \overline{B}(0,(\kappa/2+1) r)^c$,   it holds that
$$
    h(A_r)\,\ge\, c\, r^{-d} \int_{B(0,r)} h(z)\, dz\, .
$$

\end{lemma}
\pf Since $h$ is regular harmonic in
$D\cap \overline{B}(0, (\kappa/2+1) r)^c$ and
$B(A,\frac{\kappa r}{2})\subset D\cap \overline{B}(0, (\kappa/2+1) r)^c$, we have
\begin{eqnarray}
    h(A_r)=\E_{A_r}\left[h(X_{\tau_{B(A_r,\frac{\kappa r}{2})}})\right]
    \ge \int_{B(0,r)} K_{B(A_r,\frac{\kappa r}{2})} (A_r,z) h(z)\, dz\, .\label{e:ff1}
\end{eqnarray}
By \eqref{e:pke-lower} we have
 \begin{eqnarray}
    K_{B(
A_r,\frac{\kappa r}{2})} (
A_r,z)\ge c_1 j(|2(
A_r-z)|)\phi\left(\left(\frac{\kappa r}{2}\right)^{-2}\right)^{-1}\, ,\quad z\in B(0,r)\, .\label{e:ff2}
 \end{eqnarray}
Since for $z\in B(0,r)$ we have that $|A_r-z|< (\kappa^{-1}+1)r$,
by \eqref{e:doubling-condition} we have $ j(|2(A_r-z)|)\ge c_2 j(r)$ for some constant $c_2=c_2(\phi, \kappa)>0$.
Hence, combining \eqref{e:ff1}--\eqref{e:ff2} and applying Lemma \ref{l:phi-property}, we get
$$
h(A_r)\ge
\int_{B(0,r)}j(r)\phi(r^{-2})^{-1} h(z)\, dz
\ge
c_3 r^{-d} \int_{B(0,r)}h(z)\, dz
$$
which finishes the proof.
 \qed

\begin{corollary}\label{c:lemma3.2}
Let $D\subset \R^d$ be an open set which is $\kappa$-fat at infinity with
characteristics $(R,\kappa)$. There exists $c=c(d,\phi, \kappa)>0$ such that
for every $r
\ge R$ with $D\cap B(0,r)\neq \emptyset$ and every $w\in D\cap B(0,r)$ it holds that
\begin{equation}\label{e:cor3.3}
G_D(A_r,w)\ge c r^{-d} \int_{B(0,r)} G_D(z,w)\, dz\, .
\end{equation}
\end{corollary}

\pf Let $h(\cdot ):=G_D(\cdot, w)$. Then $h$ is regular harmonic in $D\cap \overline{B}(0,(\kappa/2+1) r)^c$ so the claim follows from Lemma \ref{l:lemma3.2}. \qed

\begin{lemma}\label{l:lemma3.3}
Let $D\subset \R^d$ be an open set which is $\kappa$-fat at infinity with
characteristics $(R,\kappa)$.
For $r>0$ and $n=0,1,2,\dots$,
let $B_n(r)=B(0,(\kappa/2)^{-n} r)$. There exist $c_1=c_1(d,\phi,\kappa)>0$ and
$c_2=c_2(d,\phi, \kappa)\in (0,1)$ such that for any
$r \ge R$ and any non-negative function $h$ which is regular harmonic in
$D\cap \overline{B}(0,r)^c$ and vanishes in $D^c\cap \overline{B}(0,r)^c$ we have
\begin{equation}\label{e:lemma3.3}
    \E_x\left[h(X_{\tau_{D\cap
    \overline{B_n(r)}^c}}):\, X_{\tau_{D\cap  \overline{B_n(r)}^c}}\in B(0,r)\right]\le c_1 c_2^n h(x)\, ,
\quad x\in D\cap  \overline{B_n(r)}^c,\ n=0,1,2,\dots
\end{equation}
\end{lemma}

\pf We fix $r \ge R$. For $n=0,1,2,\dots$, let $B_n=B_n(r)$, $\overline{B}_n=\overline{B_n(r)}$ and $\eta_n=(\kappa/2)^{-n} r$,  and define
$$
    h_n(x):= \E_x\left[h(X_{\tau_{D\cap \overline{B}_n^c}}):\, X_{\tau_{D\cap \overline{B}_n^c}}\in
B_0\right]\, ,\quad x\in D\cap \overline{B}_n^c\, .
$$
Then for $x\in D\cap \overline{B}_{n+1}^c$ we have
$$
    h_{n+1}(x)=\E_x\left[h(X_{\tau_{D\cap \overline{B}_n^c}}):\, \tau_{D\cap \overline{B}_{n+1}^c}=
    \tau_{D\cap \overline{B}_{n}^c},  X_{\tau_{D\cap \overline{B}_n^c}}\in
B_0\right]\le h_n(x)\, .
$$
Let $A_n=A_{\eta_n}$.
Then
\begin{eqnarray*}
    h_n(A_n)=\E_{A_n}\left[h(X_{\tau_{D\cap \overline{B}_n^c}}):\, X_{\tau_{D\cap \overline{B}_n^c}}\in
B_0\right]
    \le \E_{A_n}\left[h(X_{\tau_{\overline{B}_n^c}}):\, X_{\tau_{ \overline{B}_n^c}}\in
B_0\right]
    =\int_{
B_0} K_{\overline{B}_n^c}(A_n,z) h(z)\, dz\, .
\end{eqnarray*}
By Lemma \ref{l:poisson-kernel-estimate},
there exists $
c_1=c_1(\phi, \kappa)>0$ such that
$$
    K_{\overline{B}_n^c}(A_n,z)\le
    c_1\Big(|A_n-z|^{-d}\big(\phi(\eta_n^{-2})^{-1/2}\phi((\eta_n-|z|)^{-2})^{1/2} \big)+\eta_n^{-d}\Big)\, .
$$
For $z\in B_0$ and $n \ge 1$ we have that $|A_n-z|\asymp\eta_n$ and $\eta_n-|z|\asymp\eta_n$,
thus
$$
    K_{\overline{B}_n^c}(A_n,z)\le
    c_2 \eta_n^{-d}=
    c_2 (\kappa/2)^{n d}r^{-d}\, ,\quad z\in
B_0,  \, n=1,2,3\dots\, .
$$
Therefore, by Lemma \ref{l:lemma3.2} in the second inequality below and
Lemma \ref{l:lemma3.1} in the third,
we get that for $n=1,2,3\dots$,
\begin{eqnarray*}
    h_n(A_n)\le
    c_2 (\kappa/2)^{n d}r^{-d} \int_{
B_0} h(z)\, dz
    \le
    c_3 (\kappa/2)^{nd} h(A_0) \le
    c_4 (\kappa/2)^{\gamma n} h(A_n)\, ,
\end{eqnarray*}
where $\gamma\in (0, d)$ is the constant from Lemma \ref{l:lemma3.1}.
Now note that both $h_{n-1}$ and $h$ are regular harmonic in $D\cap \overline{B}_{n-1}^c$ and vanish on $B_{n-1}^c\cap D^c=B_{n-1}^c\setminus D\cap \overline{B}_{n-1}^c$.
Hence,
$$
    \frac{h_n(x)}{h(x)}\le \frac{h_{n-1}(x)}{h(x)}\le
    c_5 \frac{h_{n-1}(A_{n-1})}{h(A_{n-1})}\le
    c_5 C_2
    (\kappa/2)^{\gamma n}\, ,\quad x\in D\cap \overline{B}_n^c\, \quad
    n=2,3,4 \dots,
$$
where the second inequality follows from Corollary \ref{c:bhp-inf}.
The cases $n=0$ and $n=1$
 are clear by the harmonicity of $h$.
\qed

\begin{corollary}\label{c:corollary3.3}
Let $D\subset \R^d$ be an open set which is $\kappa$-fat at infinity with
characteristics $(R,\kappa)$.
For $r >0$ and $n=0,1,2,\dots$,
 let $B_n(r)=B(0,(\kappa/2)^{-n} r)$.
There exist $c_1=c_1(d,\phi,\kappa)>0$ and $c_2=c_2(d,\phi,\kappa)\in (0,1)$ such that for any
$r \ge R$ with $D\cap B(0,r)\neq \emptyset$, any $w\in D\cap B(0,r)$ and $n\ge 0$, we have
$$
\E_x\left[G_D(X_{\tau_{D\cap
\overline{B_n(r)}^c}},w):\, X_{\tau_{D\cap
\overline{B_n(r)}^c}}\in B(0,r)\right]\le c_1 c_2^n G_D(x,w)\, ,\quad x\in D\cap
\overline{B_n(r)}^c\, .
$$
\end{corollary}

\medskip
The following lemma is an analog of \cite[Lemma 16]{B} for
infinity.
The proof is essentially the same -- instead of using the balls that shrink to
a finite boundary point, we use concentric balls with larger and larger
radius (so they ``shrink at infinity''). Lemmas 13 and 14 from \cite{B}
are replaced by our Corollary \ref{c:bhp-inf} and Lemma \ref{l:lemma3.3} respectively.
Below we only indicate essential changes in the proof and refer the reader to the proof of \cite[Lemma 16]{B}.

\begin{lemma}\label{l:oscillation-reduction}
Let $D\subset \R^d$ be an open set which is $\kappa$-fat at infinity with
characteristics $(R,\kappa)$. There exist $c=c(d,\phi, \kappa)>0$ and $\nu=
\nu(d,\phi, \kappa)>0$ such that for
any $r\ge R$
and all non-negative functions
$u$ and $v$ on $\R^d$ which are regular harmonic in
$D\cap \overline{B}(0,r/2)^c$, vanish in $D^c\cap \overline{B}(0,r/2)^c$ and
satisfy $u(A_r)=v(A_r)$, there exists the limit
$$
I(u,v)=\lim_{|x| \to \infty,\,  x\in D}\ \frac{u(x)}{v(x)}\, ,
$$
and we have
\begin{equation}\label{e:speed}
\left| \frac{u(x)}{v(x)}-
I(u,v)\right| \le c \left(\frac{|x|}{r}\right)^{-\nu}\, ,\quad x\in D\cap \overline{B}(0,r)^c\, .
\end{equation}
\end{lemma}

\pf Let $r\ge R$ be fixed. Without loss of generality assume that $u(A_{r})=v(A_{r})=1$.
Let $n_0(d,\phi)\in \N$ to be chosen later, and let $
a=(\kappa/2)^{-n_0}$. For $n=0,1,2,\dots, $ define
$$
r_n=
a^n r,\ \ \overline{B}_n^c =\overline{B}(0,r_n)^c,\ \ \overline{D}_n^c =D\cap \overline{B}_n^c, \ \ \Pi_n=\overline{D}_n^c \setminus \overline{D}_{n+1}^c, \ \ \Pi_{-1}=B(0,r)\, .
$$
For $l=-1,0,1,\dots, n-1$ let
\begin{eqnarray}
u_n^l(x)&:=&\E_x\left[u(X_{\tau_{\overline{D}_n^c}}):\, X_{\tau_{\overline{D}_n^c}}\in \Pi_l\right]\, ,\quad x\in \R^d\, , \label{e:def-ukl}\\
v_n^l(x)&:=&\E_x\left[v(X_{\tau_{\overline{D}_n^c}}):\, X_{\tau_{\overline{D}_n^c}}\in \Pi_l\right]\, ,\quad x\in \R^d\, . \label{e:def-vkl}
\end{eqnarray}
Note that since $\Pi_l\subset B(0,
r_{l+1})$, it holds that
$$
u_n^l(x)\le \E_x\left[u(X_{\tau_{\overline{D}_n^c}}):\, X_{\tau_{\overline{D}_n^c}}\in B(0,
r_{l+1})\right]\, .
$$
Denote the constants $c_1$ and $c_2$ in Lemma \ref{l:lemma3.3} by $
{\wt C}$ and $\xi$ respectively. Apply Lemma \ref{l:lemma3.3} with
  $\tilde{r}=r_{l+1}$. Then $r_n = (\kappa/2)^{-n_0 n}r=(\kappa/2)^{-n_0(n-l-1)} \tilde{r}$,
 hence for $n=0,1,2,\dots $ and $x\in \overline{D}_n^c$,
$$
u_n^l(x)\le
{\wt C} (\xi^{n_0})^{n-l-1}u(x) \, ,\quad l=-1,0,1, \dots, n-2\, .
$$
Choose $n_0$ large enough so that $(1-\xi^{n_0})^{-1}\le 2$. Then since $\sum_{l=-1}^{n-2}(\xi^{n_0})^{n-l-1}=\xi^{n_0}\sum_{n=0}^{n-1}(\xi^{n_0})^n\le \xi^{n_0}(1-\xi^{n_0})^{-1}\le 2 \xi^{n_0}$, we have that for $n=1,2, \dots$ and $x\in \overline{D}_n^c$,
$$
\sum_{l=-1}^{n-2}u_n^l(x)\le 2
{\wt C} \xi^{n_0} u(x)\, .
$$
For any $\epsilon \in (0, 1)$, we can redefine $n_0(\epsilon, d, \phi)$
so that for $n=1,2, \dots$, $l=-1,0,1, \dots, n-2$,
$$
u_n^l(x)\le \epsilon^{n-1-l} u_n^{n-1}(x)\, ,\quad x\in \overline{D}_n^c\, .
$$
By symmetry we can also achieve that for $n=1,2, \dots$, $l=-1,0,1, \dots, n-2$,
$$
v_n^l(x)\le \epsilon^{n-1-l} v_n^{n-1}(x)\, ,\quad x\in \overline{D}_n^c\, .
$$

Now we claim that there exist constants $c_1=c_1(d,\phi, \kappa)>0$ and $\zeta=\zeta(d,\phi, \kappa)\in (0,1)$ such that for all $l=0,1,\dots, $
$$
\sup_{x\in \overline{D}_l^c}\frac{u(x)}{v(x)}\le (1+c_1 \zeta^l) \inf_{x\in \overline{D}_l^c}\frac{u(x)}{v(x)}\, .
$$
From now on the proof is essentially the same as the proof of \cite[Lemma 16]{B}, hence we omit it. \qed

\begin{remark}\label{r:oscillation-reduction}
{\rm
(i) Assume that $u$ and $v$ are nonnegative functions on $\R^d$ which are regular harmonic in
$D\cap \overline{B}(0,r/2)^c$ and vanish in $D^c\cap \overline{B}(0,r/2)^c$. Define $\wt{u}_r$ and $\wt{v}_r$ by
$$
\wt{u}_r(x):=\frac{u(x)}{u(A_r)}\, ,\qquad \wt{v}_r(x):=\frac{v(x)}{v(A_r)}\, .
$$
Then $\wt{u}_r$ and $\wt{v}_r$ satisfy assumptions of Lemma \ref{l:oscillation-reduction}, in particular $\wt{u}_r(A_r)=\wt{v}_r(A_r)$. Hence, there exists the limit
$$
I(\wt{u}_r,\wt{v}_r)=\lim_{|x| \to \infty,\,  x\in D}\ \frac{\wt{u}_r(x)}{\wt{v}_r(x)}\, .
$$
Therefore we can conclude that there exists the limit
$$
I(u,v,A_r)=\lim_{|x| \to \infty,\,  x\in D}\frac{u(x)}{v(x)}=\frac{u(A_r)}{v(A_r)}
I(\wt{u}_r, \wt{v}_r)\, .
$$
Suppose that $\rho\ge R$ is another radius such that $u$ and $v$ are regular harmonic in
$D\cap \overline{B}(0,\rho/2)^c$ and vanish in $D^c\cap \overline{B}(0,\rho/2)^c$. Then the same argument using $A_{\rho}$ instead of $A_r$ would give that there exists the limit
$$
I(u,v,A_{\rho})=\lim_{|x| \to \infty,\,  x\in D}\frac{u(x)}{v(x)}=\frac{u(A_{\rho})}{v(A_{\rho})}
I(\wt{u}_{\rho}, \wt{v}_{\rho})\, .
$$
This shows that the limit is independent of the point $A_r$.

(ii)
It easily follows from \eqref{e:speed} that there exist $c=c(d,\phi, \kappa)>0$ and $\nu=
\nu(d,\phi, \kappa)>0$ such that for
any $r\ge R$,
$$
\left|\frac{u(x)}{v(x)}-\frac{u(y)}{v(y)}\right| \le c \left| \frac{x-y}{r}\right|^{-\nu}\,  \quad \forall x,y \in D\cap \overline{B}(0,r)^c
$$
for all non-negative functions
$u$ and $v$ on $\R^d$ which are regular harmonic in
$D\cap \overline{B}(0,r/2)^c$, vanish in $D^c\cap \overline{B}(0,r/2)^c$ and
satisfy $u(A_r)=v(A_r)$.
}
\end{remark}
\bigskip

From now on $D$ will be an open set which is $\kappa$-fat at infinity with
characteristics $(R,\kappa)$. Fix $x_0\in D\cap \overline{B}(0,R)^c$ and recall that
$$
M_D(x,y)=\frac{G_D(x,y)}{G_D(x_0,y)}, \quad  x,y \in D\cap \overline{B}(0,R)^c.
$$
For $r> (2|x|\vee R)$, both functions $y\mapsto G_D(x,y)$ and $y\mapsto G_D(x_0,y)$ are
regular harmonic in $D\cap \overline{B}(0,r/2)^c$ and vanish on $D^c\cap \overline{B}(0,r/2)^c$.
Hence, as an immediate consequence of Lemma \ref{l:oscillation-reduction}
and Remark \ref{r:oscillation-reduction} (i)
we get the following theorem.

\begin{thm}\label{t:martin-kernel} For each $x\in D$ there exists the limit
$$
    M_D(x,\infty):=\lim_{y\in D,\ |y|\to \infty} M_D(x,y)\, .
$$
\end{thm}

Recall that $X^D$ is the process $X$ killed upon exiting $D$.
As the process $X^D$ satisfies Hypothesis (B) in \cite{KW}, $D$  has
a Martin boundary $\partial_M D$ with respect to $X$ satisfying the following properties:
\begin{description}
\item{(M1)} $D\cup \partial_M D$ is
a compact metric space (with the metric denoted by $d$);
\item{(M2)} $D$ is open and dense in $D\cup \partial_M D$,  and its relative topology coincides with its original topology;
\item{(M3)}  $M_D(x ,\, \cdot\,)$ can be uniquely extended  to $\partial_M D$ in such a way that
\begin{description}
\item{(a)}
$ M_D(x, y) $ converges to $M_D(x, w)$ as $y\to w \in \partial_M D$ in the Martin topology,
\item{(b)} for each $ w \in D\cup \partial_M D$ the function $x \to M_D(x, w)$  is excessive with respect to $X^D$,
\item{(c)} the function $(x,w) \to M_D(x, w)$ is jointly continuous on $D\times (D\cup \partial_M D)$ in the Martin topology and
\item{(d)} $M_D(\cdot,w_1)\not=M_D(\cdot, w_2)$ if $w_1 \not= w_2$ and $w_1, w_2 \in \partial_M D$.
\end{description}
\end{description}

In the remainder of the paper whenever we speak of a bounded or an unbounded sequence of points
we always mean in the Euclidean metric (and not in the Martin metric $d$).
\begin{defn}\label{finfMb}
A point $w\in \partial_M D$ is called a finite Martin boundary point if there exists a bounded sequence $(y_n)_{n\ge 1}$, $y_n\in D$,
converging to $w$ in the Martin topology. A point $w\in \partial_M D$ is called an infinite Martin boundary point if every sequence $(y_n)_{n\ge 1}$, $y_n\in D$, converging to $w$ in the Martin topology is unbounded.  The set of finite Martin boundary points is
denoted by $\partial_M^f D$, and the set of infinite Martin boundary points by $\partial_M^{\infty} D$.
\end{defn}

\begin{remark}\label{r:finite-bdry-point}
{\rm
Suppose that $w\in \partial_M^f D$ and let $(y_n)_{n\ge 1}\subset D$ be a bounded sequence converging to $w$ in the Martin topology. Then $(y_n)_{n\ge 1}$ has a subsequence $(y_{n_k})_{k\ge 1}$ converging to a point $y$ in the Euclidean topology. It cannot happen that $y\in D$, because in this case we would have that $\lim_{y_{n_k}\to y}M_D(x,y_{n_k})=M_D(x,y)$ implying by (M3)(d) that $y=w$.  Therefore, $y\in \partial D$ -- the Euclidean boundary of $D$. In particular, this shows that for every $\epsilon >0$, the sequence $(y_n)_{n\ge 1}$ (converging to $w\in \partial_M^f D$ in the Martin topology) can be chosen so that $\delta_D(y_n)<\epsilon$ for all $n\ge 1$.
}
\end{remark}

\begin{prop}\label{p:infinite-mb}
Let $D$ be an open set which is $\kappa$-fat at infinity. Then $\partial_M^{\infty} D$ consists of exactly one point.
\end{prop}
\pf Let $w\in \partial_M^{\infty} D$  and let $M_D(\cdot,w)$ be the corresponding Martin kernel. If the sequence $(y_n)_{n\ge 1} \subset D$ converges to $w$ in
the Martin topology,
then, by (M3)(a), $M_D(x,y_n)$ converge to $M_D(x,w)$. On the other hand, $|y_n|\to \infty$, thus by Theorem \ref{t:martin-kernel},
$$
\lim_{n\to \infty}M_D(x,y_n)=\lim_{|y_n|\to \infty}M_D(x,y_n)=M_D(x,\infty).
$$
Hence, for each $w\in \partial_M^{\infty} D$ it holds that $M_D(\cdot, w)=M_D(\cdot, \infty)$.
Since, by (M3)(d),  for two different Martin boundary points $w^{(1)}$ and $w^{(2)}$ it always holds that $M_D(\cdot, w^{(1)})\neq M_D(\cdot, w^{(2)})$, we conclude that the infinite part of the Martin boundary can be identified with the single point. \qed

From now on we use the notation $\partial_M^{\infty} D=\{\partial_{\infty}\}$ and, for simplicity,
we sometimes continue to write $M_D(x,\infty)$ for the more precise $M_D(x,\partial_{\infty})$.

We now briefly discuss some properties of the finite part of the Martin boundary.
Recall that $d$ denotes the Martin metric.
For $\epsilon >0$ let
\begin{equation}\label{e:definition-U_K}
K_{\epsilon}:=\left\{w\in \partial_M^f D: d(w,\partial_{\infty}) \ge \epsilon\right\}
\end{equation}
be a closed subset of $\partial_M D$. By the definition of the finite part of the Martin boundary,
for each $w\in K_{\epsilon}$ there exists a bounded sequence $(y_n^w)_{n\ge 1}\subset D$ such that
$\lim_{n\to \infty} d(y_n^w, w)=0$. Without loss of generality we may assume that
$d(y_n^w, w)<\frac{\epsilon}{2}$ for all $n\ge 1$.

\begin{lemma}\label{l:boundedness-U_K}
There exists $C_3=C_3(\epsilon)>0$ such that $|y_n^w|\le C_3$
for all $w\in K_{\epsilon}$ and all $n\ge 1$.
\end{lemma}
\pf
We first claim that for any sequence $(y_n)_{n\ge 1}$ of points in $D$, if $|y_n|\to \infty$, then
$\lim_{n\to \infty}d(y_n,\partial_{\infty})=0$, i.e., $(y_n)_{n\ge 1}$ converges to $\partial_{\infty}$
in the Martin topology.
Indeed, since $D\cup \partial_M D$ is a compact metric space, $(y_n)$ has  a convergent
subsequence $(y_{n_k})$. Let $w=\lim_{k\to \infty}y_{n_k}$ (in the Martin topology).
Then $\lim_{k\to \infty}M_D(\cdot, y_{n_k})=M_D(\cdot, w)$. On the other hand, from
Theorem \ref{t:martin-kernel}
and Proposition \ref{p:infinite-mb} we see that $\lim_{k\to \infty}M_D(\cdot, y_{n_k})=M_D(\cdot,\infty)=M_D(\cdot,\partial_{\infty})$. Therefore,
$M_D(\cdot,w)= M_D(\cdot,\partial_{\infty})$, which  implies
that $w=\partial_{\infty}$ by (M3)(d). Since this argument also holds for
any subsequence of $(y_n)_{n\ge 1}$, we conclude that $y_n\to \partial_{\infty}$ in the Martin topology.

Now suppose the lemma is not true. Then $\{y_n^w:\, w\in K_{\epsilon}, n\in \N\}$
contains a sequence $(y_{n_k}^{w_k})_{k\ge 1}$ such that $\lim_{k\to \infty}|y_{n_k}^{w_k}|= \infty$.
By the paragraph above, we have
that $\lim_{k\to \infty}d(y_{n_k}^{w_k},\partial_{\infty})= 0$. On the other hand,
$$
d(y_{n_k}^{w_k},\partial_{\infty})\ge d(w_k,\partial_{\infty})-d(y_{n_k}^{w_k}, w_k)\ge \epsilon-\frac{\epsilon}{2}=\frac{\epsilon}{2}\, .
$$
This contradiction proves the claim. \qed

Recall that an open set $D$ is called an exterior open set if $D^c$ is compact.
\begin{corollary}\label{c:exterior-isolated}
If $D$ is an exterior open set, then $\partial_{\infty}$ is an isolated point of $\partial_M D$.
Conversely, if $D$ is open and $\kappa$-fat at infinity, and $\partial_{\infty}$
is an isolated point of $\partial_M D$, then $D$ is an exterior open set.
\end{corollary}
\pf Suppose that $D$ is an exterior open set. Then $D$ is $\kappa$-fat at infinity, hence $\partial_M D=\partial_M^f D\cup \{\partial_{\infty}\}$. Since $D^c$ is compact we see that the Euclidean boundary $\partial D$ is bounded. We show that $\partial_M^f D$ is closed in the Martin topology. This will imply that $\{\partial_{\infty}\}$ is open in $\partial_M D$, hence isolated. Let $(w_n)_{n\ge 1}$ be a sequence in $\partial_M^f D$ which converges to $w\in \partial_M D$ in the Martin topology. For each $n\ge 1$, there exists a bounded sequence $(y^{w_n}_k)_{k\ge 1}$ such that $y^{w_n}_k \to w_n$ in the Martin topology. By Remark \ref{r:finite-bdry-point}, we can assume that $\delta_D(y^{w_n}_k)=\delta_{\partial D}(y^{w_n}_k)\le 1$ for all $n\ge 1$, $k\ge 1$.  Since $\partial D$ is compact, the family $\{y^{w_n}_k:\, n\ge 1, k\ge 1\}$ is bounded. Further, because $\lim_{n\to \infty}d(w_n,w)=0$ and $\lim_{k\to \infty} d(y^{w_n}_k, w_n)=0$, we can find a sequence $(y_k)_{k\ge 1}\subset \{y^{w_n}_k:\, n\ge 1, k\ge 1\}$ such that $\lim_{k\to \infty}d(y_k,w)=0$. Clearly, the sequence $(y_k)_{k\ge 1}$ is bounded proving that $w\in \partial_M^f D$.

Conversely, assume that $\partial_{\infty}$ is an isolated point of $\partial_M D$. Then there exists $\epsilon >0$ such that $K_{\epsilon}=\{w\in \partial_M^f D:\, d(w, \partial_{\infty})\ge \epsilon\}=\partial_M^f D$. Suppose that $D$ is not an exterior open set. Then both $D$ and $D^c$ are unbounded, and therefore $\partial D$ is unbounded as well. Hence, there exists $z\in \partial D$ such that $|z|\ge 3 C_3$ where $C_3=C_3(\epsilon)$ is the constant from Lemma \ref{l:boundedness-U_K}. We can find a sequence $(z_n)_{n\ge 1}\subset D$ such that $z_n\to z$ (in the Euclidean topology) and $2C_3 \le |z_n| $ for all $n\ge 1$. Since $D\cup \partial_M D$ is compact, there exist a subsequence $(z_{n_k})_{k\ge 1}$ and $w\in D\cup \partial_M D$ such that $z_{n_k}\to w$ in the Martin topology. Clearly, $w\in \partial_M D$, and since $(z_{n_k})$ is bounded, actually $w\in \partial_M^f D$. By Lemma \ref{l:boundedness-U_K}, it holds that $|z_{n_k}|\le C_3$
(for those $z_{n_k}$ for which $d(z_{n_k},w)\le \epsilon/2$). But this contradicts that $|z_{n_k}|\ge 2 C_3$.
\qed

We continue by showing that $M_D(\cdot,\partial_{\infty})$ is harmonic in $D$ with respect to $X$.

\begin{lemma}\label{l:mk-integrability}
For every
bounded open $U\subset \overline{U}\subset D$ and every $x\in D$, $M_D(X_{\tau_U},
\partial_{\infty})$ is $\P_x$-integrable.
\end{lemma}
\pf Let $(y_m)_{m\ge 1}$ be a sequence in $D\setminus \overline{U}$ such that $|y_m|\to \infty$. Then $M_D(\cdot, y_m)$ is regular harmonic in $U$. Hence, by Fatou's lemma,
\begin{eqnarray*}
  &&  \E_x[M_D(X_{\tau_U},
\partial_{\infty})]\,=\, \E_x[\lim_{m\to \infty}M_D(X_{\tau_U}, y_m)]
    \,\le \,\liminf_{m\to \infty}\E_x[M_D(X_{\tau_U}, y_m)]\\
    &&=\,\liminf_{m\to \infty} M_D(x,y_m)\,=\,M_D(x,
\partial_{\infty})<\infty\, .
\end{eqnarray*}
\qed

\begin{lemma}\label{l:mk-harmonic}
For each $x\in D$ and $\rho\in (0,\frac13 \delta_D(x)]$,
$$
    M_D(x,
\partial_{\infty})=\E_x[M_D(X_{\tau_{B(x,\rho)}},
\partial_{\infty})]\, .
$$
\end{lemma}
\pf Fix $x\in D$ and $\rho \in (0,\frac13 \delta_D(x)]$.
For $m\in \N$, let $\eta_m:=(\kappa/2)^{-m} \rho$. Let $\wt{m}\in \N$ be large enough
so that $\eta_{\wt{m}} \ge (2|x|+2\rho)\vee R$. In case $m \ge \wt{m}$,
let $A_m:=A_{\eta_m}$. Then for $m\ge \wt{m}$, $M_D(\cdot, A_m)$ is regular harmonic in
$D\setminus B(A_m, \kappa\eta_m)$ and $B(x, \rho) \subset D \setminus B(A_m, \kappa\eta_m)$, hence
\begin{equation}\label{e:mk-interchange}
    M_D(x,A_m)=\E_x[M_D(X_{\tau_{B(x,\rho)}},A_m)]\, ,\quad m\ge \wt{m}\, .
\end{equation}
From now on we assume that $m \ge \wt{m}$. To prove the statement of the
lemma it suffices to show that there exists $m_1\in \N$, $
m_1\ge \wt{m}$, such that the family $\left\{M_D(X_{\tau_{B(x,\rho)}},A_m):\, m\ge
m_1\right\}$
  is uniformly integrable with respect to $\P_x$.
This will allow us to exchange the order of the expectation and the limit when we
take the limit $m\to\infty$ in \eqref{e:mk-interchange}, thus proving the statement.

Choose $\wt{m}$ even larger so that $(\kappa/2)^{-\wt{m}}
\rho \ge 2|x_0|$, and let $m\ge \wt{m}$. Let $w\in D\cap B(0,\eta_m)$. Then $G_D(w,\cdot)$ is regular harmonic in $D\cap \overline{B}(0,\eta_m)^c$ and vanishes on $D^c\cap \overline{B}(0,\eta_m)^c$. The same is valid for $G_D(x_0,\cdot)$. Since $A_m\in D\cap \overline{B}(0, (\kappa+1)\eta_m)^c$, Corollary \ref{c:bhp-inf}
 implies that for  $w\in D\cap B(0,\eta_m)$,
\begin{equation}\label{e:mk-harmonic-1}
    M_D(w,A_m)=\frac{G_D(w,A_m)}{G_D(x_0,A_m)}\le C_2 \frac{G_D(w,y)}{G_D(x_0,y)}=C_2 M_D(w,y) \, ,\quad\textrm{for all }y\in D \cap \overline{B}(0,(\kappa+1)\eta_m)^c\, .
\end{equation}
Hence, by letting $|y|\to \infty$,
\begin{equation}\label{e:mk-harmonic-2}
    M_D(w,A_m)\le C_2 M_D(w,
\partial_{\infty})\, , \quad m\ge \wt{m}\, .
\end{equation}
Let $\epsilon >0$ be arbitrary.
By Lemma \ref{l:mk-integrability} and \eqref{e:mk-harmonic-2}, there exists $N_0>0$ such that
\begin{eqnarray}
    \lefteqn{\E_x[M_D(X_{\tau_{B(x,\rho)}},A_m):\, X_{\tau_{B(x,\rho)}}\in D\cap B(0,\eta_m), M_D(X_{\tau_{B(x,\rho)}},A_m) >N_0]}\nonumber \\
    &\le & C_2 \E_x[M_D(X_{\tau_{B(x,\rho)}},
\partial_{\infty}):\, C_2 M_D(X_{\tau_{B(x,\rho)}},
\partial_{\infty}) >N_0]
    \,\le \, C_2\, \frac{\epsilon}{2 C_2}\,=\,\frac{\epsilon}{2}\, . \label{e:mk-harmonic-3}
\end{eqnarray}
On the other hand,
\begin{eqnarray}
    \lefteqn{\E_x[M_D(X_{\tau_{B(x,\rho)}}, A_m):\, X_{\tau_{B(x,\rho)}}\in D\cap \overline{B}(0,\eta_m)^c]} \nonumber\\
    &=&\int_{D\cap \overline{B}(0,\eta_m)^c} M_D(v, A_m)K_{B(x,\rho)}(x,v)\, dv \nonumber\\
    &\le & c_1 \int_{D\cap \overline{B}(0,\eta_m)^c} M_D(v,A_m)j(|v-x|-\rho)\phi(\rho^{-2})^{-1}\, dv \, ,\label{e:mk-harmonic-4}
\end{eqnarray}
where the last inequality follows from the uniform upper estimate of the Poisson kernel
in \eqref{e:pke-upper-worse}. Choose $m_0\ge \wt{m}$ large enough  such that for $m\ge m_0$ and $v\in D\cap \overline{B}(0,\eta_m)^c$ it holds that $|v-x|-\rho\ge |v|/2$. Then $j(|v-x|-r)\le j(|v|/2)\le c_2 j(|v|)$ by \eqref{e:doubling-condition}. Hence, by treating $\phi(\rho^{-2})^{-1}$ as a constant (depending on $\rho$, but note that $\rho$ is fixed), we get that
\begin{eqnarray}\label{e:mk-harmonic-5}
    \lefteqn{\E_x[M_D(X_{\tau_{B(x,\rho)}},A_m):\, X_{\tau_{B(x,\rho)}}\in D\cap \overline{B}(0,\eta_m)^c]}\nonumber \\
    &\le & c_3 \int_{D\cap \overline{B}(0,\eta_m)^c} M_D(v,A_m)j(|v|)\, dv\, , \nonumber \\
    &=& c_3 G_D(x_0,A_m)^{-1}\int_{D\cap \overline{B}(0,\eta_m)^c} G_D(v,A_m)j(|v|)\, dv\, ,\quad m\ge m_0\, . \label{e:mk-harmonic-6}
\end{eqnarray}
By Lemma \ref{l:lemma3.1} (applied to $r=\eta_{m_0}$) we have that
\begin{equation}\label{e:mk-harmonic-7}
    G_D(x_0,A_m)^{-1}\le c_4 (\kappa/2)^{(-d+\gamma)(m-m_0)}
    G_D(x_0,A_{m_0})^{-1}\, ,\quad m\ge m_0\, ,
\end{equation}
where $\gamma\in (0,d)$. Now we estimate the integral in \eqref{e:mk-harmonic-6}:
\begin{eqnarray*}
    \lefteqn{\int_{D\cap \overline{B}(0,\eta_m)^c} G_D(v,A_m)j(|v|)\, dv}\\
    &\le & \int_{D\cap B(A_m,(\kappa+1/2) \eta_m)}G
    (v,A_m)j(|v|)\, dv + \int_{D\cap B(A_m,(\kappa+1/2) \eta_m)^c \cap \overline{B}(0,\eta_m)^c } G_D(v,A_m)j(|v|)\, dv\\
    &=:& I_1+I_2\, .
\end{eqnarray*}
To estimate $I_1$, note that if $v\in B(A_m,(\kappa+1/2) \eta_m)$, then $|v|\ge |A_m|-|v-A_m|\ge (\kappa+1)\eta_m-(\kappa+1/2) \eta_m=
(1/2) \eta_m =(1/2) (\kappa/2)^{-m} \rho$,
 hence $j(|v|)\le c_5 j((\kappa/2)^{-m} \rho)$ by \eqref{e:doubling-condition}. Therefore, by Theorem \ref{t:J-G}
 and \eqref{e:ie-1},
\begin{eqnarray*}
     I_1&\le &
     c_6  j((\kappa/2)^{-m} \rho)\int_{B(
     A_m, (\kappa/2+1) \eta_m)}\frac{1}{|v-A_m|^d
    \phi(|v-A_m|^{-2})}\, dv\\
    & \le &
    c_7 j((\kappa/2)^{-m} \rho) \int_0^{(\kappa/2+1) \eta_m}\frac{1}{s\phi(s^{-2})}\, ds\\
    &\le &
    c_8 j((\kappa/2)^{-m} \rho)\phi(((\kappa/2)^{-m} \rho)^{-2})^{-1}
 \,\le \,
 c_9 ((\kappa/2)^{-m} \rho)^{-d}\, .
\end{eqnarray*}
In order to estimate $I_2$, let $v\in
D\cap B(A_m,(\kappa+1/2) \eta_m)^c \cap \overline{B}(0,\eta_m)^c$. If $|v|\ge \kappa^{-1}(1-\kappa)^{-1}\eta_m$, then $|v-A_m|\ge |v|-|A_m|\ge  \kappa |v|$. If $\eta_m\le |v|< \kappa^{-1}(1-\kappa)^{-1}\eta_m$, then $|v-A_m|\ge (\kappa+1/2) \eta_m \ge  \kappa (1-\kappa)(\kappa+1/2)  |v|$.
 Thus, in any case, $G_D(v, A_m)\le
 g(|v-A_m|)\le c_{10}
 g(|v|)$ by \eqref{e:doubling-conditionG}. Therefore, by Theorem \ref{t:J-G},
\begin{eqnarray*}
      I_2 &\le &
      c_{10}\int_{D\cap B(A_m,(\kappa+1/2) \eta_m)^c \cap \overline{B}(0,\eta_m)^c }
    g(|v|)j(|v|)\, dv\\
    &\le &
    c_{10} \int_{\overline{B}(0,\eta_m)^c }
    g(|v|)j(|v|)\, dv
    \,\le \,
    c_{11}\int_{\eta_m}^{\infty}\frac{1}{s^{d+1}}\, ds =
    c_{12}\eta_m^{-d}=
    c_{12}((\kappa/2)^{-m} \rho)^{-d}\, .
\end{eqnarray*}
Hence,
\begin{equation}\label{e:mk-harmonic-8}
    \int_{D\cap \overline{B}(0,\eta_m)^c} G_D(v,
    A_m)j(|v|)\, dv\le
    c_{13}((\kappa/2)^{-m} \rho)^{-d}\, , \quad m\ge m_0\, .
\end{equation}
By combining \eqref{e:mk-harmonic-4}--\eqref{e:mk-harmonic-8} we get that
\begin{eqnarray}
    \lefteqn{\E_x[M_D(X_{\tau_{B(x,\rho)}},A_m):\, X_{\tau_{B(x,\rho)}}\in D\cap \overline{B}(0,\eta_m)^c]  }\nonumber \\
    &\le&
    c_{14} (\kappa/2)^{(-d+\gamma)(m-m_0)} G_D(x_0, A_{m_0})^{-1} ((\kappa/2)^{-m} \rho)^{-d} \nonumber\\
    &\le &
    c_{15} G_D(x_0, A_{m_0})^{-1} (\kappa/2)^{d m_0} (\kappa/2)^{\gamma(m-m_0)}\, , \label{e:mk-harmonic-9}
\end{eqnarray}
where in the last line we treat $\rho$ as a constant. Since $\gamma >0$, we can choose $m_1=m_1(\epsilon,m_0,d,\phi,\rho)> m_0$ large enough so that the right-hand side in \eqref{e:mk-harmonic-9} is less than $\epsilon/2$ for all $m\ge m_1$.

Finally, for $m\ge m_1$ we have
\begin{eqnarray*}
    \lefteqn{\E_x[M_D(X_{\tau_{B(x,\rho)}},A_m):\, M_D(X_{\tau_{B(x,\rho)}},A_m)>N_0]} \\
    &\le & \E_x[M_D(X_{\tau_{B(x,\rho)}},A_m):\, X_{\tau_{B(x,\rho)}}\in D\cap \overline{B}(0,\eta_m)^c] \\
    & & + \ \E_x[M_D(X_{\tau_{B(x,\rho)}},A_m):\, X_{\tau_{B(x,\rho)}}\in D\cap B(0,\eta_m), M_D(X_{\tau_{B(x,\rho)}},A_m) >N_0] \\
    & \le &\frac{\epsilon}{2}+\frac{\epsilon}{2} =\epsilon\, .
\end{eqnarray*}
Hence $\left\{M_D(X_{\tau_{B(x,\rho)}},A_m):\, m\ge m_1\right\}$ is
uniformly integrable with respect to $\P_x$.
\qed

\begin{thm}\label{t:mk-harmonic}
The function $M_D(\cdot,
\partial_{\infty})$ is harmonic in $D$ with respect to $X$.
\end{thm}
\pf
The proof of the theorem is exactly the same as that of \cite[Theorem 3.9]{KSV6}. \qed

\medskip
Let $x\in D$ and choose
$r\ge (2|x|\vee |x_0|)$.
By Corollary \ref{c:bhp-inf} we have that for all $y\in D\cap \overline{B}(0,r)^c$
$$
\frac{G_D(x,y)}{G_D(x_0,y)}\le
C_2 \frac{G_D(x,A_r)}{G_D(x_0,A_r)}\, .
$$
By letting $|y|\to \infty$ we get that
$$
M_D(x,
\partial_{\infty})\le
C_2 \frac{G_D(x,A_r)}{G_D(x_0,A_r)}\, .
$$
Suppose that $z\in \partial D$ is a regular boundary point.
Then $\lim_{x\to z}G_D(x,A_r)=0$ implying  also that
\begin{equation}\label{e:mk-limit-at-finite}
    \lim_{x\to z} M_D(x,
\partial_{\infty})=0\, , \quad
\textrm{ for  every regular boundary point } z\in \partial D.
\end{equation}

\begin{lemma}\label{l:irregular}
Suppose that $u$ is a bounded nonnegative harmonic function for $X^D$. If there exists a polar set $N\subset \partial D$ such that for any $z\in \partial D\setminus N$
$$
\lim_{D\ni x\to z} u(x)=0
$$
and
$$
\lim_{x\in D, |x|\to \infty} u(x)=0\, ,
$$
then $u$ is identically equal to zero.
\end{lemma}
\pf Take an increasing sequence of bounded open sets $\{D_n\}_{n\ge 1}$ satisfying $\overline{D_n}\subset D_{n+1}$ and $\cup_{n=1}^{\infty} D_n=D$. Then $\lim_{n\to \infty}\tau_{D_n}=\tau_D$ and by the quasi-left continuity, $\lim_{n\to \infty}X_{\tau_{D_n}}=X_{\tau_D}$ if $\tau_D<\infty$, and $\lim_{n\to \infty}|X_{\tau_{D_n}}|=\infty$ if $\tau_D=\infty$.

Since $N$ is polar, we have $\P_x(X_{\tau_D}\in N, \tau_D<
\infty)=0, x\in D$. By  harmonicity we have for every $x\in D$ and all $n\ge 1$
\begin{eqnarray*}
u(x)&=&\E_x\left[u(X_{\tau_{D_n}})\right]\\
&=& \E_x\left[u(X_{\tau_{D_n}}), \tau_D=\infty\right] +\E_x\left[u(X_{\tau_{D_n}}), \tau_{D_m}=\tau_D \textrm{ for some }m\ge 1\right]\\
& & +\E_x\left[u(X_{\tau_{D_n}}), \tau_{D_m}<\tau_D<\infty \textrm{ for all }m\ge 1\right]\, .
\end{eqnarray*}
By using bounded convergence theorem we get that
$$
\lim_{n\to \infty}\E_x\left[u(X_{\tau_{D_n}}), \tau_D=\infty\right]=\E_x\left[\lim_{n\to \infty}u(X_{\tau_{D_n}}), \tau_D=\infty\right]=0\, ,
$$
since $|X_{\tau_{D_n}}|\to \infty$ on $\{\tau_D =\infty\}$. Next, since $u=0$ on $D^c$,
\begin{eqnarray*}
\lefteqn{\lim_{n\to \infty}\E_x\left[u(X_{\tau_{D_n}}), \tau_{D_m}=\tau_D \textrm{ for some }m\ge 1\right]}\\
&=&\E_x\left[\lim_{n\to \infty}u(X_{\tau_{D_n}}), \tau_{D_m}=\tau_D \textrm{ for some }m\ge 1\right]\\
&=&\E_x\left[u(X_{\tau_D}), \tau_{D_m}=\tau_D \textrm{ for some }m\ge 1\right]=0\, .
\end{eqnarray*}
Finally, if $\tau_{D_m}<\tau_D<\infty$, then $\lim_{n\to \infty}X_{\tau_{D_n}}\in \partial D\setminus N$ $\P_x$-a.s. Hence
\begin{eqnarray*}
\lefteqn{\lim_{n\to\infty}\E_x\left[u(X_{\tau_{D_n}}), \tau_{D_m}<\tau_D<\infty \textrm{ for all }m\ge 1\right]}\\
&=&\E_x\left[\lim_{n\to \infty}u(X_{\tau_{D_n}})
{\bf 1}_{\{X_{\tau_D}\in \partial D\setminus N\}}, \tau_{D_m}<\tau_D<\infty \textrm{ for all }m\ge 1\right]=0\, .
\end{eqnarray*}
Therefore, $u(x)=0$ for every $x\in D$. \qed

Recall that a positive harmonic function $f$ for  $X^{D}$ is minimal if, whenever
$g$ is a positive harmonic function for $X^{D}$ with $g\le f$ on $D$,
one must have $f=cg$ for some constant $c$.

\bigskip

\pff {\bf of Theorem \ref{mM-inf}.}
It remains to show that $\partial_{\infty}$ is a minimal boundary point, i.e., that $M_D(\cdot, \partial_{\infty})$ is a minimal harmonic function.

Let $h$ be a positive harmonic function for  $X^{D}$
such that $h\le M_D(\cdot,
\partial_{\infty})$. By
the Martin representation in \cite{KW},
there is a measure on $\partial_M D=\partial_M^f D \cup\{\partial_{\infty}\}$ such that
$$
    h(x)=\int_{\partial_M D}M_D(x,w)\, \mu(dw)=\int_{\partial_M^f D}M_D(x,w)\, \mu(dw)+ M_D(x,
\partial_{\infty})\mu(\{\partial_{\infty}\})\, .
$$
In particular, $h(x_0)=\mu(\partial_M D)\le M(x_0,
\partial_{\infty})=1$ (because of the normalization at $x_0$). Hence, $\mu$ is a sub-probability measure.

For $\epsilon >0$,
$
K_{\epsilon}
$
is the compact subset of $\partial_M D$ defined in \eqref{e:definition-U_K}.
Define
$$
    u(x):=\int_{ K_{\epsilon} }M_D(x,w)\, \mu(dw).
$$
Then $u$ is a positive harmonic function with respect to  $X^{D}$
and bounded above as
\begin{align}
u(x)= h(x)-\mu(\{\infty\})M_D(x,
\partial_{\infty})\le \big(1-\mu(\{\infty\})\big)M_D(x,
\partial_{\infty})\, . \label{e:newm1}
\end{align}
We claim that $\lim_{|x|\to \infty} u(x)=0$. By Lemma \ref{l:boundedness-U_K} there exists $
C_3=
C_3(\epsilon)>0$ such that for each $w\in K_{\epsilon}$ there exists a sequence $(y_n^w)_{n\ge 1}\subset D$
converging to $w$ in the Martin topology and satisfying $|y_n^w|\le
C_3$. Without loss of generality we may assume that $
C_3\ge R$. Fix a point $x_1\in D\cap \overline{B}(0, 2
C_3)^c$ and choose an arbitrary point $y_0\in D\cap B(0,
C_3)$. Then for any $x\in D\cap \overline{B}(0,2
C_3)^c$ and any $y\in D\cap B(0,
C_3)$ we have that
\begin{eqnarray*}
\frac{G_D(x,y)}{G_D(x_0,y)}=\frac{G_D(x,y)}{G_D(x_1,y)}\, \frac{G_D(x_1,y)}{G_D(x_0,y)}
\le c_1 \frac{G_D(x,y_0)}{G_D(x_1,y_0)}\, \frac{G_D(x_1,y)}{G_D(x_0,y)}
\le  c_1 \frac{G(x,y_0)}{G_D(x_1,y_0)}\, \frac{G_D(x_1,y)}{G_D(x_0,y)},
\end{eqnarray*}
where the first inequality follows from the boundary Harnack principle,
 Theorem \ref{t:ubhp}.
Therefore for each $w\in
K_{\epsilon}$  we have
\begin{eqnarray*}
M_D(x,w)&=&\lim_{n\to \infty}\frac{G_D(x,y_n^w)}{G_D(x_0,y_n^w)}\le c_1 \frac{G(x,y_0)}{G_D(x_1,y_0)}\, \lim_{n\to \infty} \frac{G_D(x_1,y_n^w)}{G_D(x_0,y_n^w)}\\
&=& c_1\frac{G(x,y_0)}{G_D(x_1,y_0)}\ M_D(x_1, w) \le c_1\frac{G(x,y_0)}{G_D(x_1,y_0)}\ \sup_{w\in K_{\epsilon}} M_D(x_1, w)= c_2 G(x,y_0)
\end{eqnarray*}
by continuity of the Martin kernel (M3)(c).
 Now we let $|x|\to \infty$ and use that $G(x,y_0)\to 0$ to conclude that
 $\lim_{|x|\to \infty, x\in D}M_D(x,w)=0$ uniformly for
$w\in K_{\epsilon}$. By Theorem \ref{t:martin-kernel},
this and \eqref{e:newm1} immediately imply that $\lim_{|x|\to \infty} u(x)=0$.

From \eqref{e:mk-limit-at-finite} we see that $\lim_{x\to z} u(x)=0$ for every
regular $z\in \partial D$. Since the set of irregular boundary points is polar
(cf. \cite[(VI.4.6), (VI.4.10)]{BG}),
Lemma \ref{l:irregular} implies that $u\equiv 0$. This means that $\nu=\mu_{|
K_{\epsilon}}=0$. Since
$\epsilon >0$ was arbitrary and $\partial_M^f D=\cup_{\epsilon >0} K_{\epsilon}$,
we see that $\mu_{|\partial_M^f D}=0$. Hence $h=\mu(\{\partial_{\infty}\})M(\cdot,
\partial_{\infty})$ showing that $M(\cdot,
\partial_{\infty})$ is minimal. Therefore we have proved Theorem \ref{mM-inf}
\qed

At the end we briefly discuss the Martin boundary of the half-space $\bH=\{x=(\wt{x}, x_d):\, \wt{x}\in \R^{d-1}, x_d>0\}$.
Let $
V(r)$ be the renewal function of the ladder height process of one-dimensional subordinate Brownian motion $X^d_t=W^d(S_t)$.
It is known that the function $w(x):=
V((x_d)^+)$ is harmonic in $\bH$ with respect to $X$ (see \cite{KSV3}). Moreover, for every $z\in \partial \bH:=\{x=(\wt{x}, x_d):\, \wt{x}\in \R^{d-1}, x_d=0\}$ it holds that $\lim_{x\to z}w(x)=0$. Therefore we can conclude that $w$ is proportional to the minimal harmonic function $M_{\bH}(\cdot, \infty)$. In the next corollary we compute the full Martin boundary of $\bH$.

\begin{corollary}\label{c:H-inf}
The Martin boundary and the minimal Martin boundary of the half space $\bH$ with respect to $X$ can be identified with $\partial \bH \cup \{\infty\}$ and
$M_{\bH}(x,\infty)=w(x)/w(x_0)$ for $x \in \bH$.
\end{corollary}
\pf
By Theorem \ref{mM-inf} and the argument before the statement of this corollary,
we only need to show that the finite part $\partial_M^f \bH$ of the Martin boundary of $\bH$ can be identified with the Euclidean boundary $\partial \bH$ and that all points are minimal.
This was shown in \cite[Theorem 3.13]{KSV6} under the assumption that $\phi$ is comparable to the regularly varying function at infinity.
Even though this assumption is stronger
than {\bf (H1)}, using results in this paper and \cite{KSV8} (instead of using properties of regularly varying function) one can  follow the same proof line by line and show that under the
assumption {\bf (H1)} and $d >2\delta_2$, the finite part of Martin boundary
$\partial_M^f \bH$ can be identified with the Euclidean boundary $\partial \bH$ and that all points are minimal. We omit the details.
\qed

%%%%%%%%%%%%%%%%%%%%%%%%%%%%%%%%%%%%%%%%%%%%%%%%%%%%%%%%%%%%%%%%%%%%%%%%%%%%%%%%%%%%%%%%%%%%%%%%%%%%%%%%%%%%%%%%%%%%%%%%%%%%%%%%%%%%%

\vspace{.1in}
\begin{singlespace}

\end{singlespace}
\end{doublespace}

\vskip 0.1truein

{\bf Panki Kim}

Department of Mathematical Sciences and Research Institute of Mathematics,

Seoul National University, Building 27, 1 Gwanak-ro, Gwanak-gu Seoul 151-747, Republic of Korea

E-mail: \texttt{pkim@snu.ac.kr}

\bigskip

{\bf Renming Song}

Department of Mathematics, University of Illinois, Urbana, IL 61801,
USA

E-mail: \texttt{rsong@math.uiuc.edu}

\bigskip

{\bf Zoran Vondra\v{c}ek}

Department of Mathematics, University of Zagreb, Zagreb, Croatia

Email: \texttt{vondra@math.hr}


\small
\begin{thebibliography}{99}


\bibitem{BG} R. M. Blumenthal and R. K. Getoor:
{\it  Markov Processes and Potential Theory}, Academic Press,
New York, 1968.

\bibitem{B} K. Bogdan:  The boundary Harnack principle for the
fractional Laplacian.
 {\em Studia Math.}  {\bf 123(1)}(1997), 43--80.

\bibitem{BBKRSV} K. Bogdan, T. Byczkowski, T. Kulczycki,
M. Ryznar, R. Song, and Z. Vondra\v{c}ek: {\em Potential analysis
of stable processes and its extensions}.  Lecture Notes in
Mathematics, {\bf 1980}. Springer-Verlag, Berlin, 2009.

\bibitem{BKK} K. Bogdan, T. Kulczycki and M. Kwa\'snicki:
Estimates and structure of $\alpha$-harmonic functions. {\it Probab. Th. Rel. Fields}, {\bf 140}  (2008), 345--381.

\bibitem{BKuK} K. Bogdan, T. Kumagaii and M. Kwa\'snicki:
Boundary Harnack inequality for Markov processes with jumps. Preprint. 	arXiv:1207.3160.


\bibitem{Dyn} E. B. Dynkin: {\it Markov processes, Vol. I}.
Academic Press, New York, 1965


\bibitem{IW}
N.~Ikeda and S.~Watanabe: On some relations between the harmonic measure
  and the {L}\' evy measure for a certain class of {M}arkov processes, \emph{J.
  Math. Kyoto Univ.} \textbf{2} (1962), 79--95.



\bibitem{KM}
P.~Kim and A.~Mimica: Harnack inequalities for subordinate Brownian motions, \emph{
Electron. J. Probab. } \textbf{17} (2012), \#37.

\bibitem{KM2}
P.~Kim and A.~Mimica: Green function estimates for subordinate Brownian motions: stable and beyond,

{\it Trans. Amer. Math. Soc.}, to appear (2012)

\bibitem{KSV1} P. Kim, R. Song and Z. Vondra\v{c}ek: Boundary Harnack principle for subordinate Brownian motion. \ {\em Stoch. Proc. Appl.} {\bf 119} (2009), 1601--1631.

\bibitem{KSV3} P. Kim, R. Song and Z. Vondra\v{c}ek:
Potential theory of subordinated Brownian motions revisited.
{\it Stochastic analysis and applications to finance, essays in honour of Jia-an Yan}.
Interdisciplinary Mathematical Sciences - Vol. 13,  World Scientific, 2012, pp. 243--290.


\bibitem{KSV6} P. Kim, R. Song and Z. Vondra\v{c}ek: Minimal thinness for subordinate Brownian motion in half space.   \emph{Ann.~Inst.~Fourier} {\bf 62 (3)} (2012), 1045--1080.

\bibitem{KSV7} P. Kim, R. Song and Z. Vondra\v{c}ek: Uniform boundary Harnack principle for rotationally symmetric L\'evy processes in general open sets, {\it Sci. China Math.} {\bf 55}, (2012), 2193--2416.


\bibitem{KSV4} P. Kim, R. Song and Z. Vondra\v{c}ek:
Potential theory of  subordinate Brownian motions with Gaussian components.
{\em Stoch. Proc. Appl.}
{\bf 123(3)} (2013) 764--795.

\bibitem{KSV8} P. Kim, R. Song and Z. Vondra\v{c}ek:  Global uniform boundary Harnack principle
with explicit decay rate and its application.
Preprint, 2012.


\bibitem{KW} H.~Kunita and T.~Watanabe: Markov processes and Martin boundaries I,
    \emph{Illinois J.~Math.} {\bf 9(3)}  (1965) 485--526.

\bibitem{Kw}  M. Kwa\'snicki: Intrinsic ultracontractivity for stable semigroups on unbounded open sets. \emph{Potential Anal.} {\bf 31} (2009), 57--77.

\bibitem{Mi} P. W. Millar:
First passage distributions of processes with independent increments. {\it Ann. Probab.}, {\bf 3} (1975), 215--233.


\bibitem{Sat} K.-I. Sato:
{\em L\'evy Processes and Infinitely Divisible Distributions}. Cambridge University Press, Cambridge, 1999.

\bibitem{SSV} R. L. Schilling, R. Song and Z. Vondra{\v{c}}ek:
{\it Bernstein Functions: Theory and Applications}. de Gruyter
Studies in Mathematics 37. Berlin: Walter de Gruyter, 2010.

\bibitem{Sil} M. L. Silverstein: Classification of coharmonic
and coinvariant functions for a L\'evy process. {\em Ann.~Probab.}
{\bf 8} (1980), 539--575.

\bibitem{SW}
R. Song and J. Wu: Boundary Harnack principle for symmetric stable processes.
{\em J. Funct. Anal.} {168(2)} (1999), 403--427.

\bibitem{Sz1} P. Sztonyk:
On harmonic measure for L\'evy processes. {\em Probab. Math. Statist.}, {\bf 20} (2000), 383--390.



\end{thebibliography}
\end{document}